\documentclass[onefignum,onetabnum]{siamart220329}

\usepackage{lipsum}
\usepackage{amsfonts}
\usepackage{graphicx}
\usepackage{booktabs}
\usepackage{cite}
\usepackage{makecell}
\usepackage{mathtools}
\usepackage{epstopdf}
\usepackage{changes}
\usepackage{xurl}
\usepackage{amssymb}
\usepackage[figuresright]{rotating}
\usepackage{tikz}
\usetikzlibrary{decorations.pathreplacing}
\usepackage{dsfont}
\usepackage{bm}
\usepackage[bb=pazo]{mathalpha}
\usepackage{color,soul}
\usepackage{cellspace}
\usepackage{algorithm}
\usepackage{algpseudocode}

\usepackage{fixltx2e}

\MakeRobust{\Call}

\renewcommand{\algorithmiccomment}

\def\eg{\textit{e.g.}}
\def\ie{\textit{i.e.}}

\usetikzlibrary{shapes.geometric}

\newcommand{\threepointstar}[1]{%
  \resizebox{#1}{!}{%
    \begin{tikzpicture}
      \node[star, star points=3, star point ratio=4.0, draw, fill=black] {};
    \end{tikzpicture}%
  }%
}

\newcommand{\fourpointstar}[1]{%
  \resizebox{#1}{!}{%
    \begin{tikzpicture}
      \node[star, star points=4, star point ratio=4.0, draw, fill=black] {};
    \end{tikzpicture}%
  }%
}

\ifpdf
  \DeclareGraphicsExtensions{.eps,.pdf,.png,.jpg}
\else
  \DeclareGraphicsExtensions{.eps}
\fi

\algnewcommand{\Inputs}[1]{%
  \State \textbf{Inputs:}
  \Statex \hspace*{\algorithmicindent}\parbox[t]{.8\linewidth}{\raggedright #1}
}
\algnewcommand{\Outputs}[1]{%
  \State \textbf{Outputs:}
  \Statex \hspace*{\algorithmicindent}\parbox[t]{.8\linewidth}{\raggedright #1}
}
\algnewcommand{\Initialize}[1]{%
  \State \textbf{Initialize:}
  \Statex \hspace*{\algorithmicindent}\parbox[t]{.8\linewidth}{\raggedright #1}
}
\renewcommand{\algorithmiccomment}[1]{\hfill$\triangleright$\protect{\footnotesize{\textit{#1}}}}

\newsiamremark{remark}{Remark}
\newsiamremark{hypothesis}{Hypothesis}
\crefname{hypothesis}{Hypothesis}{Hypotheses}
\newsiamthm{claim}{Claim}

  \definecolor{masoncolor}{rgb}{0.98, 0.27, 0.62}

\headers{Dynamical processes on metric networks}{L. B\"ottcher and M. A. Porter}

\title{Dynamical processes on metric networks\thanks{Submitted to the editors DATE.
\funding{LB received funding from the Army Research Office through grant W911NF-23-1-0129.}}}

\author{Lucas B\"ottcher\thanks{Department of Computational Science and Philosophy, Frankfurt School of Finance and Management, 60322, Frankfurt am Main, Germany and Dept.\ of Medicine, University of Florida, Gainesville, FL, 32610, United States of
America
  (\email{l.boettcher@fs.de}).}
\and Mason A.\ Porter\thanks{Department of Mathematics, University of California, Los Angeles, CA, 90095, United States of
America; Department of Sociology, University of California, Los Angeles, CA, 90095, United States of
America; Sante Fe Institute, Santa Fe, NM, 87501, United States of America
  (\email{mason@math.ucla.edu}).}}

\usepackage{amsopn}

\makeatletter
\newcommand*{\addFileDependency}[1]{
  \typeout{(#1)}
  \@addtofilelist{#1}
  \IfFileExists{#1}{}{\typeout{No file #1.}}
}
\makeatother

\ifpdf
\hypersetup{
  pdftitle={Dynamical processes on metric networks},
  pdfauthor={L. B\"ottcher and M. A. Porter}
}
\fi

\begin{document}

\maketitle

\begin{abstract}
The structure of a network has a major effect on dynamical processes on that network. Many studies of the interplay between network structure and dynamics have focused on models of phenomena such as disease spread, opinion formation and changes, coupled oscillators, and random walks. In parallel to these developments, there have been many studies of wave propagation and other spatially extended processes on networks. These latter studies consider \emph{metric networks}, in which the edges are associated with real intervals. Metric networks give a mathematical framework to describe dynamical processes that include both temporal and spatial evolution of some quantity of interest --- such as the concentration of a diffusing substance or the amplitude of a wave --- by using edge-specific intervals that quantify distance information between nodes. Dynamical processes on metric networks often take the form of partial differential equations (PDEs). In this paper, we present a collection of techniques and paradigmatic linear PDEs that are useful to investigate the interplay between structure and dynamics in metric networks. We start by considering a time-independent Schr\"odinger equation. We then use both finite-difference and spectral approaches to study the Poisson, heat, and wave equations as paradigmatic examples of elliptic, parabolic, and hyperbolic PDE problems on metric networks. Our spectral approach is able to account for degenerate eigenmodes. In our numerical experiments, we consider metric networks with up to about $10^4$ nodes and about $10^4$ edges. A key contribution of our paper is to increase the accessibility of studying PDEs on metric networks. Software that implements our numerical approaches is available at \url{https://gitlab.com/ComputationalScience/metric-networks}.
\end{abstract}




\section{Introduction}
\label{sec:intro}
The study of dynamical processes on networks has led to many insights into the interplay between structure and dynamics~\cite{porter2016dynamical,masuda2017random,newman2018networks,d2019explosive,golubitsky2023dynamics}. For example, in models of disease spread~\cite{pastor2015epidemic,bottcher2017critical}, opinion 
dynamics~\cite{flache2017models,tian2023dynamics}, and coupled oscillators~\cite{kuramoto1975self,rodrigues2016kuramoto}, researchers have derived conditions for bifurcations and phase transitions between qualitatively different behaviors.
These results have often been accompanied by insights into the effectiveness of specific interventions and how various types of failures affect the robustness of structures such as communication networks~\cite{albert2000error,cohen2000resilience,cohen2001breakdown,granell2013dynamical,masoomy2023impact}, heterogeneous materials~\cite{bonamy2011failure,berthier2019forecasting,pournajar2022edge}, and social networks~\cite{danon2012social,danon2013social,schneider2022epidemic,xia2022controlling}. In these applications, the dynamical processes on the networks often take the form of
ordinary differential equations (ODEs), with each node of a network associated with one or more ODEs, which describe how their states evolve. Alternatively, many of these dynamical processes have then the form of difference equations or stochastic processes.

In parallel to these developments, a large body of literature has focused on metric networks~\cite{kuchment2003quantum,friedman2004wave,kuchment2005quantum,berkolaiko2013introduction,andrade2016green,hofmann2021asymptotics,Berkolaiko2017,arioli2018finite,kurasov2019understanding,stoll2019optimization,porter2020nonlinearity,mehandiratta2021optimal,hofmann2021pleijel,bolin2022gaussian,blechschmidt2022comparison,brio2022spectral,kravitz2022metric,kravitz2023localized}. A \emph{metric network}\footnote{We use the term ``metric network'' instead of ``metric graph'' to strengthen the link between our work and the network-science literature, where the term ``network'' is much more common than the word ``graph'' and also sometimes refers to more general objects than ordinary graphs. Some works (see, \eg, \cite{kuchment2003quantum,kuchment2005quantum,berkolaiko2013introduction}) use the term ``quantum graph'', as they specifically consider Schr\"odinger operators on metric networks.} consists of a combinatorial graph $\mathcal{G} = (\mathcal{V},\mathcal{E})$ along with a metric, where each edge $(u,v) \in \mathcal{E}$ that connects a pair of nodes $u,v \in \mathcal{V}$ is associated with a real interval $[0, \ell_{uv}]$ of length $\ell_{uv}$. If we explicitly know the position $\mathbf{x}_u$ of each node $u \in \mathcal{V}$, then $\ell_{uv} = \|\mathbf{x}_u - \mathbf{x}_v\|$ 
for a suitable norm $\| \cdot \|$, such as the Euclidean norm or (more generally) a $p$-norm. Metric networks encompass a wide variety of networked systems in which distance information between nodes is necessary to mathematically describe corresponding physical, chemical, or biological processes. Because of the edge-specific intervals $[0,\ell_{uv}]$, one can equip a metric network with a differential operator, rather than a discrete operator (such as the combinatorial Laplacian), as in a combinatorial network. This allows one to study partial differential equations (PDEs) on networks. There are also some papers that study PDEs, such as reaction--diffusion systems, on combinatorial networks (\eg, see \cite{asllani2014turing,asllani2015turing}). There is also some research on PDEs on graphons (see, \eg, \cite{medvedev2023}), which one can obtain in limiting situations from combinatorial networks.


\subsection{Prior research on metric networks and related systems}\label{prior}

In Table~\ref{tab:examples}, we overview models and application areas that are associated with metric networks. Given the wide range of scientific domains that include {PDEs on} metric networks, we only overview a small portion of the available literature. For another summary of application areas, see Chapter 7 of \cite{berkolaiko2013introduction}.

As an illustration, consider a spring network in which each edge $(u,v)$ is a spring. The end points of $(u,v)$ are nodes with positions $x_u$ and $x_v$, which we assume for simplicity are located on a line. {This example does not yield a PDE on a network, but it enables us to (1) motivate the use of edge-specific length intervals in metric networks and (2) establish connections between metric networks and combinatorial networks (\ie, the usual type of network).} The length of the edge $(u,v)$ is $\ell_{uv} = |x_u - x_v|$, and we use $w_{uv}$ to denote the corresponding spring constant. By Hooke's law, the force that acts on these nodes $u$ and $v$ is $w_{uv}\ell_{uv}$~\cite{forrow2018functional}. We fix a subset $\mathcal{W} \subseteq \mathcal{V}$ of nodes in space and seek to determine the equilibrium positions of the remaining nodes, which are in the complementary subset $\mathcal{V} \setminus \mathcal{W}$. In equilibrium, the potential energy $\sum_{(u,v)\in \mathcal{E}} w_{uv}(x_u - x_v)^2/2 = \sum_{(u,v)\in \mathcal{E}} x_u x_v L_{uv}/2$ is minimized. The matrix $L$, which has entries $L_{uv} = -w_{uv}$ for $u \neq v$ and $L_{uu} = \sum_{v\neq u} w_{uv}$, is known as the combinatorial graph Laplacian \cite{newman2018networks}. For nodes $u \in \mathcal{V} \setminus \mathcal{W}$, one achieves this minimum when
\begin{equation}
    \sum_{(u,v) \in \mathcal{E}} w_{uv} \underbrace{(x_u - x_v)}_{= \ell_{uv}} = \sum_{(u,v)\in \mathcal{E}}x_vL_{uv} = 0\,. 
\end{equation}
The lengths $\ell_{uv}$ in our spring-system example of a metric network are thus the equilibrium distances when all spring forces balance each other.

\begin{table}
\footnotesize
\centering
\begin{tabular}{S{m{3.8cm}} S{m{5.5cm}} S{m{2cm}}}\toprule\multicolumn{1}{c}{\textbf{Model}} & \multicolumn{1}{c}{\textbf{Comments}} & \multicolumn{1}{c}{\textbf{References}} \\ \hline
\makecell[l]{Quantum graphs}  & \makecell[l]{Models of quantum dyanmics \\ in thin structures.} &  \cite{ruedenberg1953free,alexander1983superconductivity,kottos1997quantum,kottos1999periodic,kuchment2002graph,kuchment2003quantum,kuchment2005quantum,berkolaiko2013introduction,Berkolaiko2017} \\ \hline
\makecell[l]{Transmission line and \\ electrical networks}  & \makecell[l]{Such networks consist of resistance, \\inductance, capacitance, and conductance\\ elements.} &  \cite{paul2007analysis,strub2019modeling,alonso2017power,chen2017power,muranova2019notion,muranova2020eigenvalues,muranova2021effective,muranova2022effective,muranova2022networks} \\ \hline
Traffic flow on networks & \makecell[l]{Flow models of vehicular and pedestrian \\ traffic, telecommunication networks, and \\ supply chains.} & \cite{piccoli2006traffic,d2010modeling,helbing2013verkehrsdynamik} \\ \hline
Gas networks  & \makecell[l]{Distribution networks that consist \\ of pipes, valves, compressors, and \\ heating and cooling elements.} &  \cite{banda2006gas,brouwer2011gas,domschke2011combination,mindt2019entropy,DomschkeHillerLangetal.2021} \\ \bottomrule
\end{tabular}
\vspace{1mm}
\caption{A variety of models and application areas of metric networks.}
\label{tab:examples}
\end{table}

Spring networks are common in models of engineering and material structures~\cite{hrennikoff1941solution,ashurst1976microscopic,beale1988elastic,hassold1989brittle,herrmann1989fracture,gusev2004finite,bottcher2021computational} and in
computer graphics~\cite{DBLP:journals/cgf/NealenMKBC06,DBLP:series/synthesis/2018Stuyck}. Networks of masses and springs
have also inspired the development of both the Gaussian-network model~\cite{bahar1997direct,haliloglu1997gaussian,cui2005normal} and the anisotropic-network model~\cite{doruker2000dynamics,atilgan2001anisotropy}, which are used to model macromolecules.

The original focus in research on metric networks concentrated on Schr\"odinger equations on networks.
The linear Schr\"odinger equation plays a central role in studies of quantum graphs, in which one uses metric graphs and considers quantum dynamics in thin structures~\cite{ruedenberg1953free,alexander1983superconductivity,kottos1997quantum,kottos1999periodic,kuchment2002graph,kuchment2003quantum,kuchment2005quantum,berkolaiko2013introduction,Berkolaiko2017}.
One can also incorporate a cubic nonlinearity to obtain a cubic nonlinear Schr\"odinger (NLS) equation, which is paradigmatic system with diverse applications. 
It arises via a mean-field description of Bose--Einstein condensates~\cite{pitaevskii2016bose}, as an envelope equation in optics~\cite{malomed2006encyclopedia}, and in many other situations. 
In the context of metric networks, cubic NLS equations have been considered on a Y-junction~\cite{noja2014nonlinear}, a dumbbell network~\cite{marzuola2016ground}, and star networks~\cite{kairzhan2019drift}. Other studies of nonlinear PDEs on metric networks include examinations of a nonlinear Dirac equation (a relativistic wave equation) on a Y-junction~\cite{sabirov2018dynamics}, the sine--Gordon equation on star and tree networks \cite{sobirov2016}, the Korteweg--de Vries equation on a star graph \cite{cavalcante2018}, and reaction--diffusion equations \cite{wallace2017}.

Metric networks also arise in the mathematical description of transmission-line and electrical networks. Such networks are usually described by lumped-element models with resistance, inductance, capacitance, and conductance elements arranged in a network, through which signals can propagate~\cite{paul2007analysis,strub2019modeling,alonso2017power,chen2017power,muranova2019notion,muranova2020eigenvalues,muranova2021effective,muranova2022effective,muranova2022networks}. For example, in Sections 22.6 and 22.7 of \cite{feynman1971feynman}, Feynman used an infinite ladder network that consists of capacitors and inductors to illustrate the function of a low-pass filter that prevents the propagation of high-frequency modes of an electromagnetic wave. In linear transmission lines, the propagation of electromagnetic waves is described by the telegraph equation~\cite{paul2007analysis}.
Because of the mathematical similarities between the telegraph equation and Schr\"odinger systems on metric networks, quantum graphs have been studied experimentally using transmission-line networks~\cite{hul2004experimental,lawniczak2019non}.

Metric networks have also appeared prominently in other applications.  
For example, networks of resistors have been used to model composite materials that consist of a combination of conducting and nonconducting materials~\cite{kirkpatrick1973percolation,zabolitzky1984monte,herrmann1984superconductivity,de1985anomalous}. In the context of quantum communication networks~\cite{nokkala2023}, 
information is transmitted through metric networks, such as optical fibers.
Metric networks are also commonly used in transport processes, including the flow of traffic, supplies, and gas in infrastructure and distribution networks~\cite{piccoli2006traffic,d2010modeling,helbing2013verkehrsdynamik,banda2006gas,brouwer2011gas,domschke2011combination,mones2014shock,mindt2019entropy,DomschkeHillerLangetal.2021}.

There are also related dynamical processes that, while not described by PDEs, are spatially extended and often arise through discretizations of PDEs. Examples of such dynamical systems include nonlinear lattice systems~\cite{chirikov1979universal, kaloshin2014arnol, flach2008discrete, kartashov2011,kotwal2021} and models of networked oscillators that have {been} used to construct classical analogs of topological insulators~\cite{susstrunk2015observation, salerno2016floquet} and spin--orbit coupling~\cite{salerno2017spin}. The Ablowitz--Ladik model~\cite{ablowitz1975nonlinear, ablowitz1976nonlinear, xia2019ablowitz} is a network of nonlinear elements that arises by discretizing an NLS equation. Other nonlinear lattice models, which are also relevant to study on more general network structures, include Fermi--Pasta--Ulam--Tsingou (FPUT) lattices~\cite{fermi1955studies, Grant_LANL, flaschka1974toda, flaschka1974todaII, dauxois2005fermi, li2023recurrence} and Toda lattices~\cite{toda1967vibration, toda1976development, teschl2000jacobi, toda2012theory}.


\subsection{Our contributions}\label{contributions}
  
  The study of PDEs on metric networks has focused on very small networks thus far~\cite{porter2020nonlinearity}. A key reason is that it is challenging to develop robust numerical methods to solve different types of PDE problems and account for diverse boundary conditions on such networks. Discretization of PDEs can yield large systems of equations that are difficult to handle numerically, especially for metric networks with many edges and PDEs that require a very small step size. Alternatively, one can employ spectral methods, although it is also challenging to identify characteristic wavenumbers and eigenmodes with high numerical precision.
  
In the present paper, we study Poisson, heat, and wave equations as paradigmatic examples of linear elliptic, parabolic, and hyperbolic
 PDE problems on metric networks. Building on previous work~\cite{arioli2018finite,gaio2019nanophotonic,brio2022spectral,kravitz2023localized}, we present different simulation approaches that are useful to study such linear PDEs on metric networks. Specifically, we extend the spectral approach of~\cite{gaio2019nanophotonic,brio2022spectral} to account for degenerate eigenmodes. 
Complementing the numerical results by Brio et al.~\cite{brio2022spectral}, who examined the Poisson equation and the telegraph equation on a metric network with three nodes, we study the Poisson equation, heat equation, and wave equation on three separate metric networks. Our numerical computations {use sparse-matrix} representations, which allow us to study metric networks with up to about $10^4$ nodes and about $10^4$ edges. A key contribution of our paper is to increase the accessibility of studying PDEs on metric networks.


\subsection{Organization of our paper}\label{organization}

Our paper proceeds as follows. In Section~\ref{sec:metric_networks}, we define metric networks and discuss common boundary conditions in the study of PDEs on metric networks. We also present an illustrative example with a Schr\"odinger operator in a two-node network. In Section~\ref{sec:methods}, we overview numerical and analytical
methods that are useful to study
metric networks. In Section~\ref{sec:star_graph}, we examine the Schr\"odinger equation on a star network as an illustrative example.
In Section~\ref{sec:numerical_example}, we study Poisson, heat, and wave equations on metric networks. In Section~\ref{sec:conclusion}, we summarize and discuss our results.
 In Appendix~\ref{app:symmetries}, we discuss group-theoretic methods that can help identify degenerate eigenmodes in metric networks.\footnote{In such degenerate eigenmodes, the corresponding wavenumbers have an algebraic multiplicity that is larger than 1.}
 In Appendix~\ref{app:larger_scale}, we solve the Poisson equation on metric networks with up to about $10^4$ nodes and about $10^4$ edges.
Our code for our numerical simulations is available at \url{https://gitlab.com/ComputationalScience/metric-networks}. 
%


\section{Metric networks}
\label{sec:metric_networks}
In Section~\ref{sec:basic_definitions}, we present some basic definitions. In Sections~\ref{sec:operators} and \ref{sec:boundary_conditions}, we overview different operators and boundary conditions in the study of PDEs on metric networks. As an illustrative example, in Section~\ref{sec:two-node}, we consider a Schr\"odinger problem in a two-node system. We point out a connection between certain solutions of this problem and a particle in an infinite square well. 
%


\subsection{Basic definitions}
\label{sec:basic_definitions}

\begin{figure}
    \centering
    \includegraphics[width=0.7\textwidth]{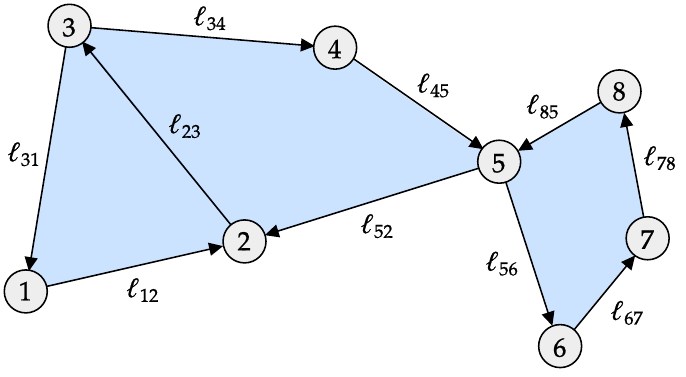}
    \caption{An example of a metric network with $N = 8$ nodes and $M = 10$ edges. The length of an edge that connects two nodes $u$ and $v$ is $\ell_{uv}$. 
        The arrows indicate the starting and ending points of an interval. For example, the edge that connects nodes 1 and 2 starts at node 1 and ends at node 2.
    The depicted metric network is embedded in $\mathbb{R}^2$, so each node is associated with a location in the plane. The blue regions indicate polygons whose vertices correspond to the network nodes. The edges of these polygons correspond to the intervals (but without their directions).
    }
    \label{fig:network_example}
\end{figure}

Consider a network in the form of a graph $\mathcal{G} = (\mathcal{V},\mathcal{E})$, where $\mathcal{V}$ is a set of nodes and $\mathcal{E}$ is a set of edges. We use $N = |\mathcal{V}|$ and $M = |\mathcal{E}|$ to denote the numbers of nodes and edges, respectively. In a \emph{metric network}, each edge $(u,v)\in \mathcal{E}$ that connects two nodes $u,v\in \mathcal{V}$ is parameterized by an interval $[0,\ell_{uv}]$ such that $0 < \ell_{uv} < \infty$. (Some authors also account for infinite-length edges~\cite{kuchment2003quantum}.) This turns the combinatorial graph $\mathcal{G}$ (which, in other contexts, is often called simply a ``graph" \cite{newman2018networks}) into a topological and metric structure. We allow loops (\ie, self-edges) and multiple edges (\ie, multi-edges) between the same nodes. The length of a 
walk that is associated with edges $\mathcal{E}' \subseteq \mathcal{E}$ is $\sum_{(u,v)\in \mathcal{E}'} \ell_{uv}$. 
For example, for the metric network in Fig.~\ref{fig:network_example}, the length of the walk $1 \rightarrow 2 \rightarrow 3 \rightarrow 4$ is $\ell_{12} + \ell_{23} + \ell_{34}$. 
In addition to intervals $[0,\ell_{uv}]$, one can equip edges with weights $w_{uv}$, as is the case in spring networks, where $w_{uv}$ represents the spring stiffness (see Section~\ref{prior}). 
In some applications, it is useful to consider time-dependent edge lengths $\ell_{uv}(t)$.

Metric networks do not have to be embedded in Euclidean space. However, several of the applications in Section~\ref{prior} and Table~\ref{tab:examples} naturally require such an embedding (\eg, gas networks, transmission lines, and quantum dynamics in thin structures) in real physical implementations of them. One can interpret a metric network as a one-dimensional (1D) simplicial complex that consists of 1D simplices (\ie, edges). A key difference between simplicial complexes in combinatorial networks~\cite{bick2023higher} and those in metric networks is that the latter are geometric entities that include length information.

Because of the interval information in a metric network, it has not only discrete nodes connected by edges but also includes all intermediate points between those edges. This allows us to define an $L^2$ space $L^2(\mathcal{G})\coloneqq \bigoplus_{(u,v)\in\mathcal{E}}L^2((u,v))$ that is associated with a metric network $\mathcal{G}$. 
Each edge $(u,v)$ has an associated continuous function $f_{(u,v)}\colon [0,\ell_{uv}] \rightarrow \mathbb{R}$ 
that maps $x \in [0,\ell_{uv}]$ to the real numbers. We require these functions to be square integrable for all edges. 
That is,
\begin{equation}
    \|f\|^2_{L^2(\mathcal{G})} \coloneqq \sum_{(u,v) \in \mathcal{E}}\|f_{(u,v)}\|^2_{L^2\left((u,v)\right)} < \infty\,,
\end{equation}
where
\begin{equation}
    \|f_{(u,v)}\|^2_{L^2\left((u,v)\right)} \coloneqq \langle f_{(u,v)},f_{(u,v)}\rangle_{L^2\left((u,v)\right)} = \int_0^{\ell_{uv}} f_{(u,v)}(x)^2 \,\mathrm{d}x\,.
\end{equation}
We calculate inner products between two functions $f,g \in L^2(\mathcal{G})$ by computing
\begin{equation}
 \langle f,g\rangle_{L^2\left(\mathcal{G}\right)}\coloneqq \sum_{(u,v)\in \mathcal{E}} \langle f_{(u,v)}, g_{(u,v)}\rangle_{L^2\left((u,v)\right)}=\sum_{(u,v)\in \mathcal{E}} \int_0^{\ell_{uv}} f_{(u,v)}(x) g_{(u,v)}(x)\,\mathrm{d}x\,.
\end{equation}
%


\subsection{Operators}
\label{sec:operators}
Because of the edge-specific intervals $[0,\ell_{uv}]$, one can equip a metric network with differential operators rather than discrete operators (such as the combinatorial Laplacian), which are studied often in combinatorial networks.

Relevant operators $\mathcal{H}\colon H^2([0,\ell_{uv}])\rightarrow H^2([0,\ell_{uv}])$ that act on $f_{(u,v)}$ include the negative second derivative
\begin{equation}
    \mathcal{H}(f_{(u,v)})(x) = -\frac{\mathrm{d^2}}{\mathrm{d}x^2}f_{(u,v)}(x)\,,
    \label{eq:d2fdx2}
\end{equation}
the Schr\"odinger operator
\begin{equation}
    \mathcal{H}(f_{(u,v)})(x) = -\frac{\mathrm{d^2}}{\mathrm{d}x^2}f_{(u,v)}(x)+U(x)f_{(u,v)}(x)\,,
\end{equation}
and the magnetic Schr\"odinger operator~\cite{auletta2012}
\begin{equation}
    \mathcal{H}(f_{(u,v)})(x)=\left(-\mathrm{i}\frac{\mathrm{d}}{\mathrm{d}x} - A(x)\right)^2 f_{(u,v)}(x)+U(x)f_{(u,v)}(x)\,,
\end{equation}
where $U\colon [0,\ell_{uv}]\rightarrow \mathbb{R}$ is a scalar potential function and $A \colon [0,\ell_{uv}]\rightarrow \mathbb{R}$ is a vector potential function. 
The space $H^2([0,\ell_{uv}])$ is the Sobolev space of twice-differentiable functions on the interval $[0,\ell_{uv}]$. Sobolev spaces arise commonly in the analysis (including numerical analysis) of PDE problems in their weak formulations \cite{evans2010}.
For a metric network with edges $(u,v) \in \mathcal{E}$, the function space is $\bigoplus_{(u,v)\in \mathcal{E}} H^2([0,\ell_{uv}])$.
%


\subsection{Boundary conditions}
\label{sec:boundary_conditions}
To solve a PDE on a metric network, we need to incorporate suitable boundary conditions for all functions $f_{(u,v)}(x)$ at their end points (\ie, for $x\in\{0,\ell_{uv}\}$). We first require that $f$ {satisfies a continuity condition} on $\mathcal{G}$. That is, for each node $u$ with degree $\mathrm{deg}(u)$, it needs to satisfy
${\mathrm{deg}(u) - 1}$ equations that ensure the
continuity of all $\mathrm{deg}(u)$ functions $f_{(u,v)}(x)$.
Additionally, for each node $u$, it is common to employ the Kirchhoff flux condition
\begin{equation}
    \sum_{e \in \mathcal{E}_u}\left.\frac{\mathrm{d} f_{e}(x)}{\mathrm{d}x} \right|_{x = x_e^*} = 0\,,
    \label{eq:kirchhoff}
\end{equation}
where $\mathcal{E}_u$ denotes the set of edges that are attached to node $u$ and we choose $x_e^*$, where $e \in \mathcal{E}_u$ is an edge, to evaluate $f_e(\cdot)$ at node $u$. As in~\cite{kuchment2003quantum,Berkolaiko2017},
we use the convention that derivatives are taken away from a node into an edge. Some works refer to the Kirchhoff condition as the ``Kirchhoff--Neumann'' condition or the ``Neumann'' condition. (See, \eg, \cite{kuchment2003quantum} and \cite{Berkolaiko2017}.) The connection to the Neumann condition in standard PDE problems becomes apparent if we consider a node $u$ with a single incident edge $e$. In this case, Eq.~\eqref{eq:kirchhoff} requires that the derivative of $f_e(\cdot)$ vanishes at the node $u$.

Let $\widetilde{\mathcal{H}}$ denote the operator that includes both the relevant derivatives and the boundary conditions on $\mathcal{G}$. In the present paper, we focus on problems that involve negative second derivatives of $f_{(u,v)}(x)$ with respect to $x$ [see Eq.~\eqref{eq:d2fdx2}] and primarily consider the Kirchhoff flux condition~\eqref{eq:kirchhoff}. 
It has been shown~\cite{kostrykin1999kirchhoff,kuchment2003quantum,berkolaiko2013introduction} that the resulting operator $\widetilde{\mathcal{H}} = -\tilde{\Delta}$ is self-adjoint and hence has an orthonormal eigenbasis and real eigenvalues. This is a key result in the study of metric networks, as it allows one to expand PDE solutions in the eigenbasis of $\widetilde{\mathcal{H}}$.

Another boundary condition that preservers the self-adjointness of the Schr\"odinger operator~\eqref{eq:d2fdx2} is the Dirichlet condition
\begin{equation}     \label{eq:dirichlet}
    \left. f_e(x)\right|_{x=x_e^*} = 0\quad \text{for all}\quad e \in \mathcal{E}_u\,.
\end{equation}
Imposing Dirichlet conditions at each node yields a metric network that consists of noninteracting edges.

In the present paper, we use the term ``coupling conditions'' to refer to the combination of continuity conditions
and (either Kirchhoff or Dirichlet) boundary conditions {for all nodes}.


\subsection{Two-node system}
\label{sec:two-node}
As an illuminating example, we examine a PDE on a simple metric network. Consider the linear, time-independent Schr\"odinger equation on a network of two nodes that are connected by a single edge of length $\ell$.\footnote{Henceforth, when we use the term ``Schr\"odinger equation'', we always mean the linear, time-independent Schr\"odinger equation (\ie, the Helmholtz equation).} We seek to determine the solution $f(x)$ of the Schr\"odinger (\ie, Helmholtz) equation 
\begin{equation}
    \frac{\mathrm{d}^2f(x)}{\mathrm{d}x^2}=-k^2 f(x)\,,\quad x\in[0,\ell]\,.
    \label{eq:example}
\end{equation}
The boundary condition is given by the Kirchhoff flux condition~\eqref{eq:kirchhoff}, which yields $f'(0)=0$ and $-f'(\ell) = 0$. 
For completeness (and pedantry), we include the minus sign in the boundary condition at $x = \ell$, following the convention that derivatives are taken away from a node into an edge. The solution of Eq.~\eqref{eq:example} is $f(x) = A e^{\mathrm{i}kx} + Be^{-\mathrm{i}kx}$. The boundary condition $f'(0) = 0$ implies that $A = B$. Because $-f'(\ell) = 0$, we obtain $k_m = \pi m/\ell$ with $m \in\{1,2,\dots\}$. In quantum mechanics, 
$m$ is known as a ``quantum number''; in this
example, $m$ labels the vibrational modes of a particle in a box. We discard the trivial solution $f(x) \equiv 0$. The eigenfunctions that are associated with the eigenvalues $k_m$ are 
\begin{equation}    \label{eq:example_kirchhoff_wall}
    f(x;m) = 2A\cos(\pi m x/\ell)\,.
\end{equation}
Adding a node between the two existing nodes in the interval $[0,\ell]$ does not change the solution \eqref{eq:example_kirchhoff_wall}. A degree-two node with the Kirchhoff flux condition is thus equivalent to an uninterrupted edge~\cite{kuchment2003quantum,Berkolaiko2017}.

\begin{figure}
    \centering
    \includegraphics{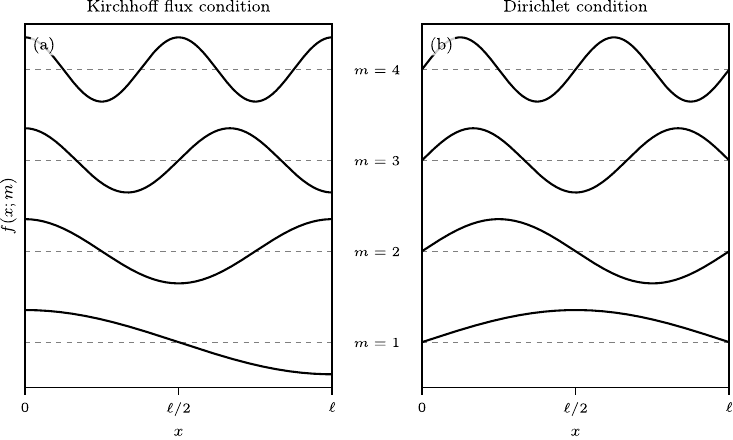}
    \caption{The eigenfunctions associated with Eq.~\eqref{eq:d2fdx2} on a metric network with two nodes that are connected by a single length-$\ell$ edge.
     The eigenfunctions are (a) $f(x;m) = \sqrt{2/\ell}\cos(k_m x)$ for Kirchhoff flux boundary conditions and (b) $f(x;m) = \sqrt{2/\ell}\sin(k_m x)$ for Dirichlet boundary conditions, where $k_m = (\pi m)/\ell$ is the wavenumber for $m \in \{1, 2, \ldots\}$.
     }
    \label{fig:eigenfunctions_line}
\end{figure}

One can interpret the function $f(x;m) = 2A\cos(\pi m x/\ell)$ as a quantum
wave function [see Fig.~\ref{fig:eigenfunctions_line}(a)]. This requires us to normalize it such that
\begin{equation}
    \int_0^\ell |f(x;m)|^2\,\mathrm{d}x = 1\,,
\end{equation}
which yields $A = 1/\sqrt{2 \ell}$ and
\begin{equation}     \label{eq:example_kirchhoff_wall_norm}
    f(x;m) = \sqrt{\frac{2}{\ell}}\cos(\pi m x/\ell)\,.
\end{equation}

As an alternative to the Kirchhoff flux condition \eqref{eq:kirchhoff}, one can employ the Dirichlet condition \eqref{eq:dirichlet}. A calculation that is similar to the one above yields
\begin{equation}     \label{eq:example_infinite_well}
    f(x;m) = \sqrt{\frac{2}{\ell}}\sin(\pi m x/\ell)\,.
\end{equation}
Equation~\eqref{eq:example_infinite_well} describes a quantum particle in an infinite square well [see Fig.~\ref{fig:eigenfunctions_line}(b)].

Consistent with the self-adjointness of the negative second derivative with Kirchhoff flux conditions (see Section~\ref{sec:boundary_conditions}), both Eqs.~\eqref{eq:example_kirchhoff_wall_norm} and \eqref{eq:example_infinite_well} give orthonormal bases with respect to the inner product $\langle f_{(u,v)}, g_{(u,v)}\rangle_{L^2\left((u,v)\right)}$ (see Section~\ref{sec:basic_definitions}). That is, 
\begin{equation}
    \int_0^{\ell} f(x;m) f(x;m')\,\mathrm{d}x = \delta_{mm'}\,,    
\end{equation}
where the Kronecker delta function $\delta_{mm'}$ is 1 for $m = m'$ and is $0$ otherwise.
%


\section{Numerical and analytical methods}
\label{sec:methods}
There are different approaches to solve PDEs on metric networks. It is sometimes possible to obtain analytical solutions for small networks and PDEs that are analytically tractable on each edge. For larger networks and/or analytically intractable PDEs, it is necessary to employ numerical methods to obtain solutions. In Sections \ref{sec:sm} and \ref{sec:fd}, we discuss two numerical approaches: a spectral method~\cite{gaio2019nanophotonic,brio2022spectral} and a finite-difference method. Using a spectral method to solve a PDE on a metric network with a self-adjoint operator (see our illustrative example in Section~\ref{sec:operators}) seems especially suitable, given the availability of a Fourier-like basis. Other numerical methods to study PDEs on metric networks include a discontinuous Galerkin method~\cite{brio2022spectral} and a finite-element method~\cite{arioli2018finite}. In Section~\ref{sec:weyls_law}, we discuss Weyl's law~\cite{weyls_law,Berkolaiko2017} as a way to help track wavenumbers when employing a spectral approach. Weyl's law gives an analytical estimate of the number of eigenmodes of the Schr\"odinger equation on a metric network up to a specified cutoff value. 
In Appendix~\ref{app:symmetries}, we discuss group-theoretic methods that can help identify degenerate eigenmodes in metric networks. Characterizing potential degeneracies of eigenmodes is useful when applying spectral approaches and Weyl's law.


\subsection{A spectral method}
\label{sec:sm}
One can express the solution of a linear PDE on a metric network in terms of an expansion (a so-called ``spectral expansion") with respect to
an appropriate orthonormal basis. 
One can construct such a basis using the eigenmodes of the self-adjoint operator $-\tilde{\Delta}$, the generalized negative Laplacian operator that includes continuity conditions and boundary conditions at all nodes. 

To compute the eigenmodes and corresponding wavenumbers, we solve an eigenvalue problem that accompanies the Schr\"odinger equation
\begin{equation}     \label{eq:helmholtz}
    \tilde{\Delta} f = -k^2 f\,.
\end{equation}
The function 
\begin{equation}
    f \coloneqq \left(f_1(x_1),f_2(x_2),\ldots,f_M(x_M)\right)^\top
\end{equation}
includes all functions $f_i(x_i)$ that are defined on their associated
edges, which have domains $[0,\ell_i]$. To make our notation more concise,
we write $f_i(x_i)$ (with $i\in\{1,\ldots,M\}$) instead of $f_{(u,v)}(x_{(u,v)})$ (with $(u,v)\in \mathcal{E}$) and write $\ell_i$  instead of $\ell_{(u,v)}$.
Some works (see, \eg, \cite{brio2022spectral,kravitz2023localized}) use the same argument $x$ for different edges, but we employ the notation $x_i$ (with $i \in \{1,\ldots,M\}$) to account for the possibility that different edges can
have distinct domains.

Solving Eq.~\eqref{eq:helmholtz} yields
\begin{equation}
    f_i(x_i) = A_i \sin(k x_i) + B_i \cos(k x_i)\,,
\end{equation}
where one determines the coefficients $A_i$ and $B_i$ using the imposed continuity conditions and boundary conditions (i.e., the coupling conditions). For a node $u$ with degree $\mathrm{deg}(u)$, there are $\mathrm{deg}(u) - 1$ equations 
associated with the continuity condition and $1$ equation  associated with the boundary condition.
The total number of equations is thus $\sum_{u\in \mathcal{V}}\mathrm{deg}(u) = 2M$. These equations yield the homogeneous system 
\begin{equation}     \label{eq:homogenous_system}
    T(k)X = 0
\end{equation}
for the coefficient vector $X = (A_1,B_1,\ldots,A_M,B_M)^\top\in \mathbb{R}^{2M}$. The nontrivial solutions of Eq.~\eqref{eq:homogenous_system} require the coupling-condition matrix $T(k)$ to be singular (\ie, $\mathrm{det}(T(k)) = 0$). 
We refer to the corresponding values $k_m$ (with quantum number $m \in \{1,2,\ldots\}$) as the ``characteristic wavenumbers'' of the metric network.\footnote{We prefer the term ``characteristic wavenumber'' to ``resonant frequency''~\cite{brio2022spectral,kravitz2023localized} because the quantity $k$ is physically a wavenumber, rather than a frequency.}

For each characteristic wavenumber $k_m$, we determine the nullspace of $T(k_m)$. If its dimension $\dim(\ker T(k_m))$ is larger than $1$, there exist degenerate eigenmodes. We denote the corresponding coefficients by $A^{mn}_{1},B^{mn}_1,\ldots,A^{mn}_M,B^{mn}_M$ (with $n \in \{1,\dots,\dim(\ker T(k_m))\}$). The eigenmodes that are associated with $k_m$ are
\begin{equation}
    f^{mn} =
    \begin{pmatrix}
    A^{mn}_{1}\sin(k_m x_1)+B^{mn}_{1}\cos(k_m x_1) \\
    A^{mn}_{2}\sin(k_m x_2)+B^{mn}_{2}\cos(k_m x_2) \\
    \vdots \\
    A^{mn}_{M}\sin(k_m x_M)+B^{mn}_{M}\cos(k_m x_M)
    \end{pmatrix}\,.
\end{equation}
One can normalize the eigenmodes $f^{mn}$ using the inner product
\begin{equation}
 \langle f^{mn}, f^{mn}\rangle_{L^2\left(\mathcal{G}\right)} = \sum_{i=1}^M \int_0^{\ell_i} f^{mn}_{i}(x) f^{mn}_{i}(x)\,\mathrm{d}x\,.
\end{equation}

Given a set $\{f^{mn}\}$ of orthonormal eigenmodes, we can construct a spectral expansion for another function that is defined on the same metric network. 
It has been shown~\cite{brio2022spectral} that the spectral-expansion coefficients decay faster than any polynomial (a property that is known as ``spectral convergence'') for functions in $L^2(\mathcal{G})$ with compact support on all edges for which this function is nonzero. 
For functions in $L^2(\mathcal{G})$ that do not have compact support on all such edges,
the expansion coefficients decay with the quantum number $m$ as $m^{-4}$.

In our numerical calculations, we use the described spectral approach to construct the solutions of PDEs on a metric network using a spectral expansion in the eigenmodes $\{f^{mn}\}$.  We give further details in Section~\ref{sec:poisson}, where we consider the Poisson equation on several metric networks.
%


\subsection{Finite differences}
\label{sec:fd}
\begin{figure}
    \centering
    \begin{tikzpicture}
      \def\li{8} 
      \def\hi{1} 
      
      \draw (0,0) node[left] {0} -- (\li,0) node[right] {$\ell_i$};
      
      \foreach \x/\label in {0/x_{i,0}, 1/x_{i,1}, 2/x_{i,2}, 7/x_{i,N_i-1}, 8/x_{i,N_i}} {
        \draw (\x*\hi,0.1) -- (\x*\hi,-0.1) node[below] {$\label$};
      }
    
      \node[below] at (4.5*\hi, -0.1) {$\ldots$};
      \draw[decorate,decoration={brace,amplitude=5pt}] (1*\hi, 0.2) -- (2*\hi, 0.2) node[midway, above=4pt] {$h_i\coloneqq \frac{\ell_i}{N_i}$};
    \end{tikzpicture}
    \caption{A finite-difference discretization of edge $i$ with length $\ell_i$. We employ a uniform discretization with step size $h_i = \ell_i/N_i$.}
    \label{fig:fd_discretization}
\end{figure}
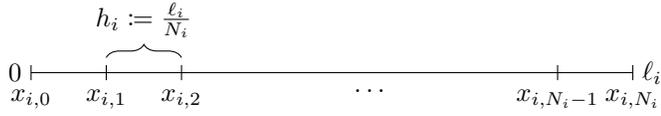

Finite-difference approximations are a complementary method to solve a PDE on a metric network.
 In the problems that we study in the present paper, we need to discretize both first derivatives (because of Kirchhoff boundary conditions) and second derivatives (\eg, for Schr\"odinger operators). We denote the step size in a discretized edge domain $[0, \ell_i]$ by $h_i \coloneqq \ell_i/N_i$ (with $i \in \{1,\dots,M\}$), where $N_i$ is the number of intervals that we use to discretize $[0, \ell_i]$. 
 In Fig.~\ref{fig:fd_discretization}, we show a schematic illustration of our discretization scheme.

One possible discretization of the second derivative of $f_i(x_i)$ is
\begin{equation}     \label{eq:d2f_dx2_fd}
    \frac{\mathrm{d}^2}{\mathrm{d}x^2}f_i(x_i)\equiv f_i''(x_i) = \frac{f_i(x_i + h_i) - 2f_i(x_i) + f_i(x_i - h_i)}{h_i^2} + \mathcal{O}(h_i^2)\,.
\end{equation}
Naturally, one can also employ higher-order discretizations or use implicit methods.
As a shorthand notation for $f_i(jh_i)$ (with $j \in \{0,\ldots,N_i\}$), we write $f_{i,j}$.

Consider a Schr\"odinger equation of the form \eqref{eq:example} on each of the $M$ edges of a metric network. We summarize the discretized second derivatives \eqref{eq:d2f_dx2_fd} in a square matrix. Most of the matrix elements are $0$, so we use a sparse matrix representation when implementing our numerical solvers. 

For Dirichlet boundary conditions, $f_i(0) = f_{i,0} = 0$ and $f_i(\ell_i) = f_{i,N_i} = 0$ for all edges $i$. 
Therefore, the second derivatives at $x_i = h_i,\ell_i-h_i$ are $f_{i,1}'' = (f_{i,2} - 2f_{i,1})/h_i^2 + \mathcal{O}(h_i^2)$ and $f_{i,N_i-1}'' = (f_{i,N_i - 2} - 2f_{i,N_i -1})/h_i^2 + \mathcal{O}(h_i^2)$. For one edge with Dirichlet boundaries, the discretized version of the generalized Laplacian $\tilde{\Delta}$ [see Eq.~\eqref{eq:helmholtz}] is thus
\begin{equation}     \label{eq:discrete_laplace_dirichlet}
    \tilde{\Delta}_{h_i} = \frac{1}{h_i^2}
    \begin{pmatrix}
  -2 & 1 & 0 & \cdots & 0 & 0 \\
  1 & -2 & 1 & \cdots & 0 & 0 \\
  0 & 1 & -2 & \cdots & 0 & 0 \\
  \vdots & \vdots & \vdots & \ddots & \vdots & \vdots \\
  0 & 0 & 0 & \cdots & 1 & -2
\end{pmatrix}\in\mathbb{R}^{(N_i-1) \times (N_i-1)}\,.
\end{equation}
The eigenvector that is associated with the discretized Schr\"odinger equation \eqref{eq:helmholtz} and discrete Laplace--Dirichlet operator \eqref{eq:discrete_laplace_dirichlet} is 
\begin{equation}
    f_{i,j'} \propto \sin\left(\frac{\pi m j'}{N_i}\right)\,, \quad j' \in \{1,\ldots,N_i - 1\}\,,
\end{equation}
which is a discrete analogue of the sine eigenfunction~\eqref{eq:example_infinite_well}. The corresponding eigenvalues $k_{i,m}$ satisfy
\begin{equation}
    -k_{i,m}^2 = -\frac{4}{h_i^2} \sin\left(\frac{\pi m}{2 N_i}\right)^2=-\frac{m^2 \pi^2}{\ell_i^2}+\mathcal{O}(m^4 h_i^2)\,, \quad m \in \{1,\ldots,N_i - 1\}\,.
    \label{eq:discrete_laplace_eigenvalue}
\end{equation}
These eigenvalues yield the eigenvalues for the continuous problem \eqref{eq:example} in the limit $h_i \rightarrow 0$~\cite{chung2000discrete,brio2022spectral}.

For a single edge with Kirchhoff (\ie, Neumann) boundaries, 
$f'_{i,0} = (f_{i,1} - f_{i,-1})/(2 h_i) + \mathcal{O}(h_i^2) = 0$ and $f'_{i,N_i} = (f_{i,N_i - 1} - f_{i,N_i + 1})/(2h_i) + \mathcal{O}(h_i^2) = 0$. We have introduced ``ghost'' points at the two positions $x_i = -h_i$ and $x_i = \ell_i + h_i$ to write a second-order finite-difference approximation at the boundaries.\footnote{As emphasized in \cite{brio2022spectral}, it is important to maintain uniform approximation orders both in an edge and at its boundaries.
The lowest-order approximation determines the overall order of the employed approximation scheme.} The discretized version of the generalized Laplacian $\tilde{\Delta}$ for a single Kirchhoff edge is
\begin{equation}
    \tilde{\Delta}_{h_i} = \frac{1}{h_i^2}
    \begin{pmatrix}
  -2 & 2 & 0 & \cdots & 0 & 0 \\
  1 & -2 & 1 & \cdots & 0 & 0 \\
  0 & 1 & -2 & \cdots & 0 & 0 \\
  \vdots & \vdots & \vdots & \ddots & \vdots & \vdots \\
  0 & 0 & 0 & \cdots & 2 & -2
\end{pmatrix}\in\mathbb{R}^{(N_i+1)\times (N_i+1)}\,.
    \label{eq:discrete_laplace_kirchhoff}
\end{equation}
Solving the discrete Schr\"odinger eigenvalue problem \eqref{eq:helmholtz} using the discrete Laplace--Kirchhoff operator \eqref{eq:discrete_laplace_kirchhoff} yields
\begin{equation}
    f_{i,j'} \propto \cos\left(\frac{\pi m j'}{N_i}\right)\,, \quad j' \in \{0,\ldots,N_i\}\,,
\end{equation}
which is a discrete analogue of the cosine eigenfunction~\eqref{eq:example_kirchhoff_wall_norm}.
The corresponding eigenvalues  $k_{i,m}$ satisfy \eqref{eq:discrete_laplace_eigenvalue} with $m \in \{0,\ldots,N_i\}$. In signal processing and data compression, the matrices \eqref{eq:discrete_laplace_dirichlet} and \eqref{eq:discrete_laplace_kirchhoff} are known as the discrete sine transform and discrete cosine transform, respectively~\cite{strang1999discrete}.

\begin{figure}
    \centering
    \includegraphics{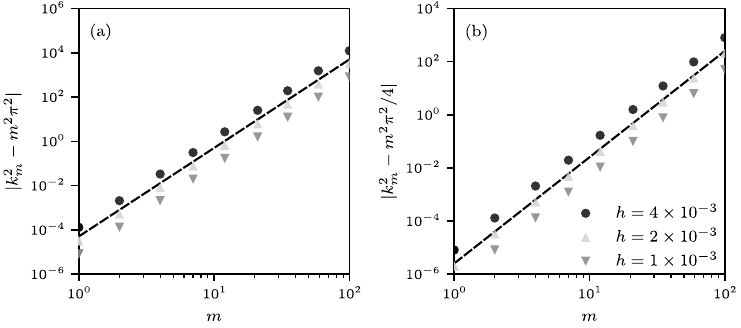}
    \caption{Error scaling of our finite-difference method for the Schr\"odinger equation on metric networks. We show the absolute error in the difference between the numerically obtained $k_{m}^2$ (where $k_m$ is the eigenvalue) and the corresponding analytical values of $k_m$ as a function of the quantum number $m$. Different types of markers correspond to different step sizes $h$. (a) The error scaling for a metric network with two nodes, a single length-$1$ edge, and Dirichlet boundaries.
     (b) The error scaling for a star network with 4 nodes, 3 length-$1$ edges, and Kirchhoff flux boundaries.
     In both examples, the observed error scaling is consistent with a quartic dependence on
     $m$, as indicated by the dashed black line [see Eq.~\eqref{eq:discrete_laplace_eigenvalue}].}
    \label{fig:fd_error_scaling}
\end{figure}

It is straightforward to simulate the Schr\"odinger equation \eqref{eq:example} on a metric network with Dirichlet boundaries. One just needs to construct a block-diagonal matrix in which each block represents the Laplace--Dirichlet operator \eqref{eq:discrete_laplace_dirichlet} that is associated with a specified edge. Recall that Dirichlet boundaries imply that edges are isolated, resulting in noninteracting ``signals''. The situation is different for metric networks with Kirchhoff flux boundaries. To construct the discretized generalized Laplacian for a metric network with Kirchhoff boundaries, one possible starting point is to use a block-diagonal matrix in which each block represents the Laplace--Kirchhoff operator \eqref{eq:discrete_laplace_dirichlet} that is associated with a specified edge. One then needs to adjust the matrix entries such that edges interact through Kirchhoff flux conditions at the associated nodes. Consider a node $u$ at which the edges $i \in \mathcal{E}_u$ terminate or originate. We use $f_0$ to denote the value of $f$ at the node $u$. Regardless of the edge's orientation, we use $f_{i,1}$ to denote the value of the function $f_i$ at the discretization point next to node $u$. As in \cite{brio2022spectral}, using a central second-order scheme to approximate the first derivative at node $u$ yields
\begin{equation}    \label{eq:kirchhoff_fd}
   2 \frac{\sum_{i\in\mathcal{E}_u} \frac{f_{i,1} - f_{0}}{h_i}}{\sum_{i\in \mathcal{E}_u} h_i} = -k^2 f_{i,0}\,.
\end{equation}
For each node with Kirchhoff boundaries, one needs to incorporate the associated expression from the left-hand side of Eq.~\eqref{eq:kirchhoff} into the generalized discretized Laplacian matrix.

Equation~\eqref{eq:kirchhoff_fd} gives one way to couple the dynamics that are associated with different edges. In a recent paper~\cite{avdonin2022discretization}, Avdonin et al.\ used a variational approach to establish coupling conditions for the discretized wave equation on a metric network.

Although finite differences provide a relatively straightforward way to numerically solve the Schr\"odinger equation \eqref{eq:helmholtz} on a metric network, a downside of this approach is the $\mathcal{O}(m^4 h_i^2)$ error term in $k_{i,m}^2$ [see Eq.~\eqref{eq:discrete_laplace_eigenvalue}] for both Dirichlet and Kirchhoff boundary conditions.

In Fig.~\ref{fig:fd_discretization}, we show the absolute error in
the difference between the numerically obtained $k_{i,m}^2$ and the corresponding analytical values as a function of $m$. We consider two metric networks: a 2-node network with a single edge and Dirichlet boundaries [see Fig.~\ref{fig:fd_error_scaling}(a)] and a star network with $N = 4$ nodes, $M = 3$ edges, and Kirchhoff flux boundaries [see Fig.~\ref{fig:fd_error_scaling}(b)]. 
The observed error scaling is consistent with the aforementioned quartic dependence on $m$. 

\begin{figure}
    \centering
    \includegraphics[width=0.39\textwidth]{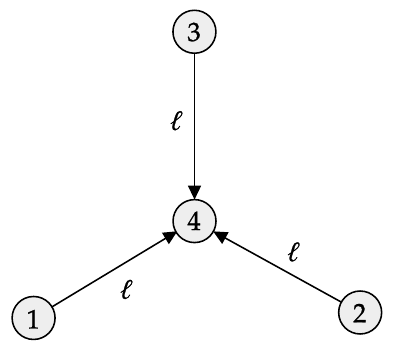}
    \caption{A metric star network with $N = 4$ nodes and $M = 3$ edges. The length of each edge is $\ell$.}
    \label{fig:star_graph}
\end{figure}

We show a schematic illustration of a 3-edge star network in Fig.~\ref{fig:star_graph}. We will revisit this example in Section~\ref{sec:star_graph} and in Appendix~\ref{app:symmetries}. In the present discussion, our objective is to emphasize that employing finite differences may not be practical when attempting to capture signals with large wavenumers (\ie, small wavelengths). However, this approach can be useful in situations with small wavenumbers. It can also provide benchmark results to use as a baseline when employing other numerical techniques (such as the spectral method in Section~\ref{sec:sm}).
%


\subsection{Weyl's law}
\label{sec:weyls_law}
Identifying all characteristic wavenumbers in a given interval can be challenging because of the rounding errors that occur when numerically determining if $T(k)$ [see Eq.~\eqref{eq:homogenous_system}] becomes singular and when working with discretized Laplace operators (see Section~\ref{sec:fd}). Therefore, it is useful to estimate of the number of characteristic wavenumbers up to a specified cutoff value.
Let $N_\mathcal{G}(k)$ denote the characteristic-wavenumber counting function. This function counts the number of characteristic wavenumbers $k'$ that satisfy $k'^2 \leq k^2$. That is,
\begin{equation}
    N_\mathcal{G}(k)\coloneqq\{|k' \colon k'^2 \leq k^2 \,\,\, \mathrm{and} \,\,\,
    \mathrm{det}(T(k')) = 0|\}\,.    
\end{equation}
According to Weyl's law~\cite{berkolaiko2013introduction,Berkolaiko2017}, 
\begin{equation}
    N_\mathcal{G}(k) = \frac{\mathcal{L}}{\pi}k + \mathcal{O}(1)\,,
\end{equation}
where $\mathcal{L} = \sum_{i=1}^M \ell_{i}$ is the total length of the edges. Additionally, the counting function $N_\mathcal{G}(k)$ satisfies 
\begin{equation}
    \frac{\mathcal{L}}{\pi}k-M\leq N_\mathcal{G}(k) \leq \frac{\mathcal{L}}{\pi}k+N\,.
    \label{eq:weyl_bounds}
\end{equation}
Deviations of the counting function $N_\mathcal{G}(k)$ from Weyl's law have been studied both theoretically~\cite{davies2010non,davies2012non} and experimentally~\cite{lawniczak2019non}.
%


\section{Illustrative example: A star network}
\label{sec:star_graph}
As an example, consider the Schr\"odinger equation \eqref{eq:helmholtz} on a metric star network $\mathcal{G}$ with $N = 4$ nodes, $M = 3$ {equal-length} edges, and Kirchhoff flux conditions at each node. 
We seek to compute the characteristic wavenumbers and their corresponding eigenmodes. This example is analytically tractable, but we also use a numerical spectral approach that we will continue to use subsequently.\footnote{One can obtain analytical insights even for metric star networks with three unequal edge lengths. Barra and Gaspard~\cite{barra2000level} derived an analytical description of the distribution of level spacings (\ie, the differences between consecutive energy levels) for the Schr\"odinger operator [see Eq.~\eqref{eq:helmholtz}] on such metric networks.}
Our comparison of analytical and numerical results for the considered star network enables us to examine the numerical-resolution requirements of the spectral 
approach in Section~\ref{sec:sm}.

The coupling-condition matrix $T_{\threepointstar{0.15cm}}(k)$ [see Eq.~\eqref{eq:homogenous_system}] that is associated with the 3-edge star network with Kirchhoff flux conditions is
\begin{equation}
    T_{\threepointstar{0.15cm}}(k) =
    \begin{pmatrix}
    1 & 0 & 0 & 0 & 0 & 0 \\
    0 & 0 & 1 & 0 & 0 & 0 \\
    0 & 0 & 0 & 0 & 1 & 0 \\
    \sin(k \ell) & \cos(k \ell) & -\sin(k \ell) & -\cos(k \ell) & 0 & 0 \\
    \sin(k \ell) & \cos(k \ell) & 0 & 0 & -\sin(k \ell) & -\cos(k \ell) \\
    -\cos(k \ell) & \sin(k \ell) & -\cos(k \ell) & \sin(k \ell) & -\cos(k \ell) & \sin(k \ell) \\
    \end{pmatrix}\,.
    \label{eq:T_k_star_graph}
\end{equation}
Alternatively, we can first establish that $A_1 = A_2 = A_3 = 0$ because nodes 1, 2, and 3 are degree-1 nodes with the Kirchhoff flux condition (see Section~\ref{sec:two-node}). Therefore, the eigenmodes are cosine functions. The remaining Kirchhoff and continuity conditions at node 4 yield
\begin{equation}
    \widetilde{T}_{\threepointstar{0.15cm}}(k) =
    \begin{pmatrix}
    \cos(k \ell) & -\cos(k \ell) & 0 \\
    \cos(k \ell) & 0 & -\cos(k \ell) \\
    \sin(k \ell) & \sin(k \ell) & \sin(k \ell)
    \end{pmatrix}\,. 
\end{equation}
The nontrivial solution of Eq.~\eqref{eq:homogenous_system} satisfies $\mathrm{det}(\widetilde{T}_{\threepointstar{0.15cm}}(k)) = 0$. This yields
\begin{equation}
    3 \cos(k\ell)^2 \sin(k\ell) = 0\,.
\end{equation}
The characteristic wavenumbers are thus $k_m = \pi m/(2 \ell)$ (with $m \in \{1,2,\ldots\}$). 
The algebraic multiplicity of $k_m$ is
\begin{equation}
    \mu_{T(k_m) }=
    \begin{cases}
        2\,,\quad\text{$m$ is odd}\\
        1\,,\quad\text{$m$ is even}  \,.
    \end{cases}
\end{equation}
For odd $m$, there are two degenerate eigenmodes: 
\begin{equation}
    f^{m1} = \frac{1}{\sqrt{\ell}}
    \begin{pmatrix}
        \cos(k_m x)\\
        0\\
        -\cos(k_m x)
    \end{pmatrix}
    \quad\text{and}\quad
    f^{m2}=\frac{1}{\sqrt{\ell}}
    \begin{pmatrix}
        \cos(k_m x)\\
        -\cos(k_m x)\\
        0
    \end{pmatrix}\,.
    \label{eq:star_degenerate}
\end{equation}
{The observed number of degenerate eigenmodes is equal to the dimension of one of the irreducible representations of the underlying symmetry group, which is the permutation group $S_3$ (see Appendix~\ref{app:symmetries}).}
For even $m$, the eigenmode is
\begin{equation}
    f^{m1} = \sqrt{\frac{2}{3\ell}}
    \begin{pmatrix}
        \cos(k_m x)\\
        \cos(k_m x)\\
        \cos(k_m x)
    \end{pmatrix}\,.
    \label{eq:star_not_degenerate}
\end{equation}
Observe that $\langle f^{m1},f^{m2}\rangle_{L^2\left(\mathcal{G}\right)} = 1/2$ if $m$ is odd and $\langle f^{m1},f^{m'1}\rangle_{L^2\left(\mathcal{G}\right)} = \delta_{mm'}$ if $m$ is even.

\begin{figure}
    \centering
    \includegraphics{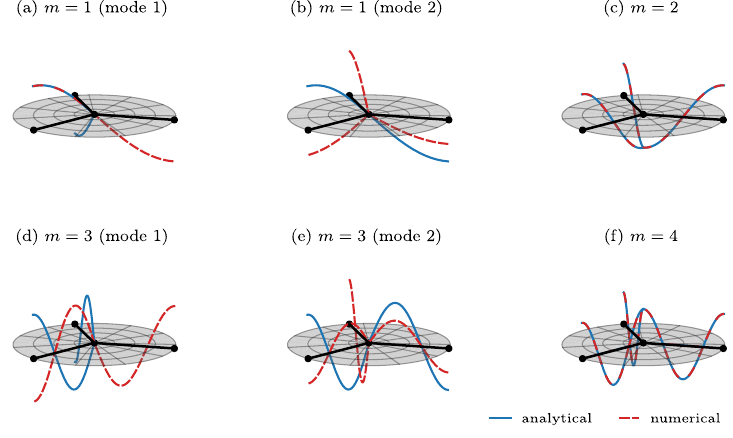}
    \caption{Eigenmodes of the Schr\"odinger equation \eqref{eq:helmholtz} on a metric star network with $N = 4$ nodes (black disks) and $M = 3$ edges (black lines). All edges have length $\ell$ and are equipped with Kirchhoff boundaries. For $m \in \{1,3\}$, the eigenmodes are degenerate [see Eq.~\eqref{eq:star_degenerate}]; they are not degenerate for $m \in \{2,4\}$ [see Eq.~\eqref{eq:star_not_degenerate}]. The solid blue curves show the analytical solutions $\eqref{eq:star_degenerate}$ and $\eqref{eq:star_not_degenerate}$. We also determine the characteristic wavenumbers $k_m$ numerically by searching for minima of the inverse condition number of the coupling-condition matrix $T_{\protect\threepointstar{0.11cm}}(k)$ [see Eq.~\eqref{eq:T_k_star_graph}]. We then use these numerical $k_m$ to compute the nullspace (and hence the eigenmodes) of the coupling-condition matrix using a QR decomposition. We indicate the associated numerical eigenmodes with dashed red curves.
    }
    \label{fig:star_network_eigenmodes}
\end{figure}

In Fig.~\ref{fig:star_network_eigenmodes}, we show the eigenmodes of the 3-edge metric star network for $m \in \{1,2,3,4\}$. In the spectral method, we seek to determine the characteristic wavenumbers $k_m$ for the entire network $\mathcal{G}$ (rather than for individual edges). Therefore, unlike in the finite-difference approach [see Eq.~\eqref{eq:discrete_laplace_eigenvalue}], we use only one index in $k_m$ when labeling an edge.
The solid blue curves show the analytically obtained eigenmodes from Eqs.~\eqref{eq:star_degenerate} and \eqref{eq:star_not_degenerate}. We also examine the ability of our spectral numerical approach to identify the characteristic wavenumbers and their corresponding eigenmodes with sufficient numerical precision. 
For each characteristic wavenumber $k_m$, the coupling-condition matrix $T(k_m)$ is singular (\ie, $\mathrm{det}(T(k_m)) = 0$). However, the determinant is not an appropriate indicator of singularity in numerical calculations.
Following \cite{gaio2019nanophotonic,brio2022spectral}, we use the condition number $\kappa(T(k))$ to determine if a certain value of $k$ is a characteristic wavenumber $k_m$. For values of $k$ that are close to $k_m$, the coupling-condition matrix is almost singular, so the condition number $\kappa(T(k))$ increases significantly as $k \rightarrow k_m$.
We compute $\kappa(T(k))$ using the identity
\begin{equation}     \label{eq:condition_number}
    \kappa(T(k)) = \frac{\sigma_{\rm max}(T(k))}{\sigma_{\rm min}(T(k))}\,,
\end{equation} 
where $\sigma_{\rm max}(T(k))$ and $\sigma_{\rm min}(T(k))$ are the maximum and minimum singular values, respectively. There exist sparse singular-value-decomposition (SVD) methods in many numerical software packages (\eg, \texttt{SciPy}) that allow one to compute the singular values of large, sparse matrices.

\begin{figure}
    \centering
    \includegraphics{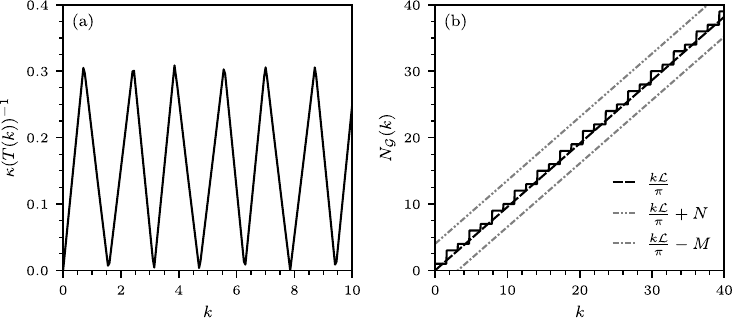}
    \caption{The inverse condition number $\kappa(T(k))^{-1}$ and characteristic-wavenumber counting function $N_\mathcal{G}(k)$ for a 3-edge metric star network. {The length of each edge is $\ell = 1$}. (a) The inverse condition number $\kappa(T(k))^{-1}$ as a function of the characteristic wavenumber $k$. Observe that $\kappa(T(k))^{-1}$ is small near the characteristic wavenumbers $k_m = \pi m/2$ (with $m \in \{1,2,\ldots\}$). (b) The characteristic-wavenumber counting function $N_\mathcal{G}(k)$ (solid black curve) as a function of $k$ and Weyl's law $N_\mathcal{G}(k) = k\mathcal{L}/\pi$ (dashed black curve), where $\mathcal{L} = 3\ell = 3$. In this plot, we include the ``zero mode" (for which $k_m = 0$) in $N_\mathcal{G}(k)$. 
    The dash-dotted and dash-dot-dotted gray lines, respectively, indicate the lower and upper bounds of Weyl's law [see Eq.~\eqref{eq:weyl_bounds}]. 
  The algebraic multiplicity of $k_m$ is 2 for odd $m$ and 
  1 for even $m$. We thus have to count $k_m$ twice for odd $m$. In numerical computations, we do this by counting $k_m$ using
  its corresponding $\dim(\ker T(k_m))$.
  }
    \label{fig:star_weyl}
\end{figure}

To numerically determine characteristic wavenumbers, we initially compute the inverse condition numbers $\kappa(T(k))^{-1}$ for a range of values of $k$ [see Fig.~\ref{fig:star_weyl}(a)]. For the star network, we consider $k \in [0,10^2]$ and choose $2 \times 10^3$ equidistant values of $k$. We consider a value of $k$ to be a candidate for a characteristic wavenumber if $\kappa(T(k))^{-1} < 10^{-2}$. We then minimize $\kappa(T(k))^{-1}$ using the candidate characteristic wavenumbers as initial values $k^{(0)}$. To minimize the scalar function $\kappa(T(k))^{-1}$, we employ a constrained Brent method~\cite{brent2013algorithms}, which is implemented in the \texttt{minimize\_scalar} function in \texttt{SciPy}. For an initial value $k^{(0)}$, we set the interval bound of $k$ to $[k^{(0)} - 0.1, k^{(0)} + 0.1]$. We use machine precision as an acceptable absolute error for the convergence of $k$.
 The best choices for
 the range and number of values of $k$, the optimization bounds, and the convergence criterion depend on the specific system that one is studying.
 After performing minimizations for all candidate characteristic wavenumbers, we obtain a set of values of $k_m$ for which $\kappa(T(k_m))^{-1}$ is close to $0$. We then round the values of $k_m$ to the nearest number with a specified number of digits and obtain a set of distinct numerical characteristic wavenumbers.

Armed with a set of characteristic wavenumbers $k_m$, we seek to compute the vectors that span the nullspace of $T(k_m)$. To do so, one can use a method that is based on a singular-value decomposition or a QR decomposition.\footnote{The sparse-matrix SVD in \texttt{SciPy} can produce erroneous nullspace vectors; see, \eg, \url{https://github.com/scipy/scipy/issues/11406}. Therefore, we use the sparse-matrix QR decomposition method at \url{https://github.com/yig/PySPQR}.} In Algorithm~\ref{alg:wavenumbers_eigenmodes}, we summarize our numerical method to determine wavenumbers and eigenmodes.
\begin{algorithm}
\footnotesize
\caption{Compute the characteristic wavenumbers and their corresponding eigenmodes.
}
\label{alg:wavenumbers_eigenmodes}
\begin{algorithmic}[1]
\Inputs{k\_arr, $\Delta k$, cutoff, precision, $T()$, $\Call{inv\_cond\_num}$, $\Call{minimize}$, $\Call{nullspace}$
}
\Outputs{k\_m\_arr, f\_m\_arr}
\Initialize{char\_wavenum\_cand, k\_m\_arr, f\_m\_arr}
\For{$k$ in k\_arr}
    \State kappa\_inv $\gets$ $\Call{inv\_cond\_num}{k,T}$ 
     \If{kappa\_inv $<$ cutoff} \quad \Comment{\,Check that the inverse condition number is sufficiently small to be a characteristic-wavenumber candidate}
        \State char\_wavenum\_cand.append($k$)
    \EndIf
\EndFor
\For{$k^{(0)}$ in char\_wavenum\_cand}
\State $k_m \gets \Call{minimize}{}(\Call{inv\_cond\_num}{},\mathrm{args} = (T),\mathrm{bounds} = (k^{(0)} - \Delta k,k^{(0)} + \Delta k))$
\State  k\_m\_arr $\gets$ $k_m$
\EndFor
\State rounded\_elements $\gets$ [\Call{round}{$k_m$, precision} \textbf{for} $k_m$ in k\_m\_arr]
\State k\_m\_arr = \Call{list}{\Call{set}{k\_m\_arr}}
\quad \Comment{\,Final set of distinct
characteristic wavenumbers}
\For{$k$ in k\_m\_arr}
\State V $\gets$ \Call{nullspace}{$k$,T}
\State f\_m\_arr.append(V)
\EndFor
\State \textbf{return} k\_m\_arr, f\_m\_arr
\end{algorithmic}
\end{algorithm}

For our 3-edge star network,  we show the numerically obtained eigenmodes as dashed red curves in Fig.~\ref{fig:star_network_eigenmodes}. 
These numerical eigenmodes coincide with the analytical eigenmodes for $m \in \{2,4\}$.
This is not the case for $m \in \{1,3\}$ because of the degeneracy of the modes. The numerical eigenmodes in Figs.~\ref{fig:star_network_eigenmodes}(a,d) are $(\cos(k_m x),-\cos(k_m x),0)^\top/\sqrt{\ell}$, which equals 
 $f^{m2}$ in  Eq.~\eqref{eq:star_degenerate}. In Figs.~\ref{fig:star_network_eigenmodes}(b,e), the numerical eigenmode
 is $(-\cos(k_m x),-\cos(k_m x),2\cos(k_m x))^\top/\sqrt{3\ell}$, which equals 
             \linebreak
 $-(2/\sqrt{3}) f^{m1} + (1/\sqrt{3}) f^{m2}$. 
 The numerically obtained degenerate eigenmodes lie in $\mathrm{span}\{f^{m1},f^{m2}\}$. 
  Unlike $f^{m1}$ and $f^{m2}$, the two numerical eigenmodes are orthonormal. 
  That is,
\begin{align}
    \left\langle -\frac{2}{\sqrt{3}}f^{m1} + \frac{1}{\sqrt{3}}f^{m2}, f^{m2} \right\rangle_{L^2\left(\mathcal{G}\right)} = &-\frac{2}{\sqrt{3}}\left\langle f^{m1}, f^{m2}\right\rangle_{L^2\left(\mathcal{G}\right)} \notag \\
    	&+ \frac{1}{\sqrt{3}}\left\langle f^{m2}, f^{m2}\right\rangle_{L^2\left(\mathcal{G}\right)} = 0\   
\end{align}
and
\begin{align}
    \left\| -\frac{2}{\sqrt{3}}f^{m1} + \frac{1}{\sqrt{3}}f^{m2} \right\|^2_{L^2\left(\mathcal{G}\right)} =  -\frac{4}{3} \left\langle f^{m1}, f^{m2} \right\rangle_{L^2\left(\mathcal{G}\right)} &+ \frac{4}{3} \left\|f^{m1}\right\|^2_{L^2\left(\mathcal{G}\right)} \notag \\
    	&+ \frac{1}{3}\left\|f^{m2}\right\|^2_{L^2\left(\mathcal{G}\right)} = 1\ . 
\end{align}
If a set of degenerate eigenmodes is not already orthonormal, one can make them orthonormal using the Gram--Schmidt algorithm.

To ensure that our numerical approach identifies
all characteristic wavenumbers in a given interval, we use Weyl's law (see Section~\ref{sec:weyls_law}) to compare the numerically obtained wavenumber-counting function $N_\mathcal{G}(k)$ to its estimate $k\mathcal{L}/\pi$. For the 3-edge star network with equal edge lengths $\ell = 1$, the total edge length $\mathcal{L}$ is $3$. In Fig.~\ref{fig:star_weyl}(b), we see that the numerically computed counting function $N_\mathcal{G}(k)$ (solid black curve) closely resembles the estimate from Weyl's law (dashed black curve). Visible differences between the two curves can highlight the need to refine a numerical method.
%


\section{Numerical examples with various PDEs}
\label{sec:numerical_example}
We now study the Poisson equation, the heat equation, and the wave equation on metric networks. These three PDEs, respectively, are fundamental types of elliptic, parabolic, and hyperbolic PDEs. They complement our study of Schr\"odinger problems in Section~\ref{sec:two-node} [see Eq.~\eqref{eq:example}] and Section~\ref{sec:methods} [see Eq.~\eqref{eq:helmholtz}].
%

\subsection{Poisson equation}
\label{sec:poisson}

\begin{figure}
    \centering
    \includegraphics{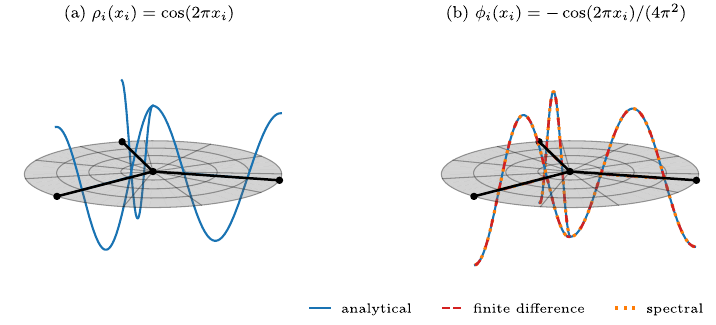}
    \caption{Solution of the Poisson equation on a metric star network with $N = 4$ nodes (black dots) and $M = 3$ edges (black lines). All edges $i \in \{1,2,3\}$ have length $\ell_i = 1$ and Kirchhoff boundaries. (a) The source term $\rho_i(x_i) = \cos(2\pi x_i)$ (solid blue curves) on the right-hand side of the Poisson equation \eqref{eq:poisson}. (b) The solution $\phi_i(x_i) = -\cos(2\pi x_i)/(4\pi^2)$ (solid blue curves) of the Poisson equation. The dashed red and dotted orange curves indicate numerical solutions using finite-difference and spectral methods, respectively.
    }
    \label{fig:poisson_star}
\end{figure}

We first consider the Poisson equation
\begin{equation}     \label{eq:poisson}
    \tilde{\Delta} \phi = \rho\,,
\end{equation}
which describes the potential field $\phi$ that is associated with a given function $\rho$ (\eg, a mass or electric-charge distribution).
In addition to its manifestations in  mechanics and electrostatics, the Poisson equation is also the steady-state equation of the heat equation with a heat source (see Section~\ref{sec:heat}).
As before, the operator $\tilde{\Delta}$ in Eq.~\eqref{eq:poisson} is the generalized Laplacian that includes continuity conditions and boundary conditions at all nodes.
The Poisson equation is an elliptic PDE.

Consider the 3-edge star network from Section~\ref{sec:methods} (see Fig.~\ref{fig:star_graph}) and set
\begin{equation}     \label{eq:source_term}
    \rho_i(x_i) = \cos(2\pi x_i)
\end{equation}
for all edges $i \in \{1,2,3\}$ [see Fig.~\ref{fig:poisson_star}(a)]. We set the lengths $\ell_i$ of all edges to $1$. All boundaries are of Kirchhoff type. The solution of the corresponding Poisson equation \eqref{eq:poisson} is $\phi_i(x_i) = -\cos(2\pi x_i)/(4\pi^2)$ [see Fig.~\ref{fig:poisson_star}(b)]. 

\begin{figure}
    \includegraphics[width=\textwidth]{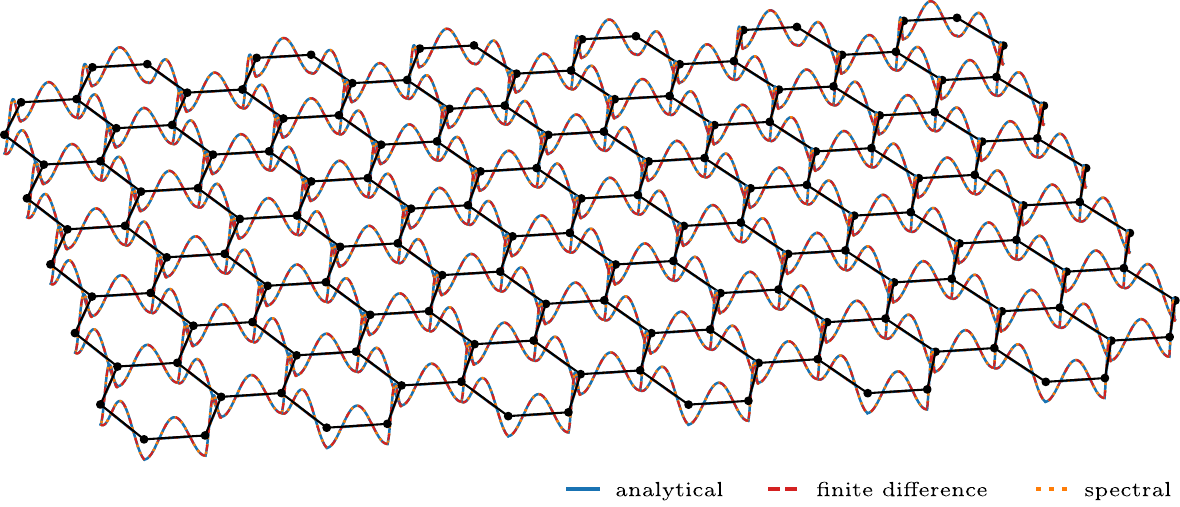}
\caption{Solution of the Poisson equation on a metric hexagonal lattice with $N = 154$ nodes (black dots) and $M = 213$ edges (black lines). All edges $i \in \{1,\ldots,M\}$ have length $\ell_i = 1$ 
and Kirchhoff boundaries. The source term on the right-hand side of the Poisson equation \eqref{eq:poisson} is $\rho_i(x_i) = \cos(2\pi x_i)$ for each edge. The solution of the Poisson equation is $\phi_i(x_i) = -\cos(2\pi x_i)/(4\pi^2)$ for each edge (solid blue curves). The dashed red and dotted orange curves indicate numerical solutions using finite-difference and spectral methods, respectively.}
\label{fig:poisson_hexagonal}
\end{figure}

Recall that the generalized negative Laplacian with Kirchhoff boundaries is self-adjoint and hence has an orthonormal eigenbasis (see Section~\ref{sec:boundary_conditions}). Therefore, we employ a spectral approach and expand the solution of Eq.~\eqref{eq:poisson} using the eigenbasis of $-\tilde{\Delta}$. To do so, we run Algorithm~\ref{alg:wavenumbers_eigenmodes} with the same parameters as in Section~\ref{sec:star_graph} to compute characteristic wavenumbers. We then construct the solution of the Poisson equation~\eqref{eq:poisson} using orthonormal spectral solutions $f^{mn}$ (with $m \in \{1,2,\ldots\}$ and $n \in \{1,\ldots,\dim(\ker T(k_m))\}$) that are associated with Eq.~\eqref{eq:helmholtz}. That is,
\begin{equation}     \label{eq:spectral_poisson}
    \phi = \sum_{m,n} a_{mn} f^{mn}\quad\mathrm{and}\quad a_{mn} = -b_{mn}/k_m^2\,,
\end{equation}
where $b_{mn} = \langle f^{mn},\rho\rangle_{L^2\left(\mathcal{G}\right)}$.\footnote{
Because the constant ``zero mode" of the eigenvalue problem~\eqref{eq:helmholtz} has an associated eigenvalue of $0$, the Fredholm alternative (as well as the compatibility condition) imply that $b_{01} = 0$~\cite{brio2022spectral}.
 Therefore, the sum over $m$ in Eq.~\eqref{eq:spectral_poisson} starts at $m = 1$.}

In Eq.~\eqref{eq:spectral_poisson}, we assume that the spectral expansion is not truncated. The summation 
over all $a_{mn} f^{mn}$ then yields the exact solution $\phi$. In practice, one has to truncate the sum at a certain value of $m$.

For our finite-difference solution of Eq.~\eqref{eq:poisson}, we set the number of discretization intervals to $N_i = 1000$ for all edges $i \in \{1,2,3\}$. To discretize the generalized Laplacian $\tilde{\Delta}$, we use the underlying Laplace--Kirchhoff matrices \eqref{eq:discrete_laplace_kirchhoff} and employ Eq.~\eqref{eq:kirchhoff_fd} to implement the Kirchhoff flux condition at the hub node at which all three edges terminate.

In Fig.~\ref{fig:poisson_star}(b), we show the numerical solutions for the 3-edge star network that we obtain using finite-difference and spectral methods. We see that both approaches are able to appropriately resolve the true solution.
\begin{figure}
    \centering
    \includegraphics[width=\textwidth]{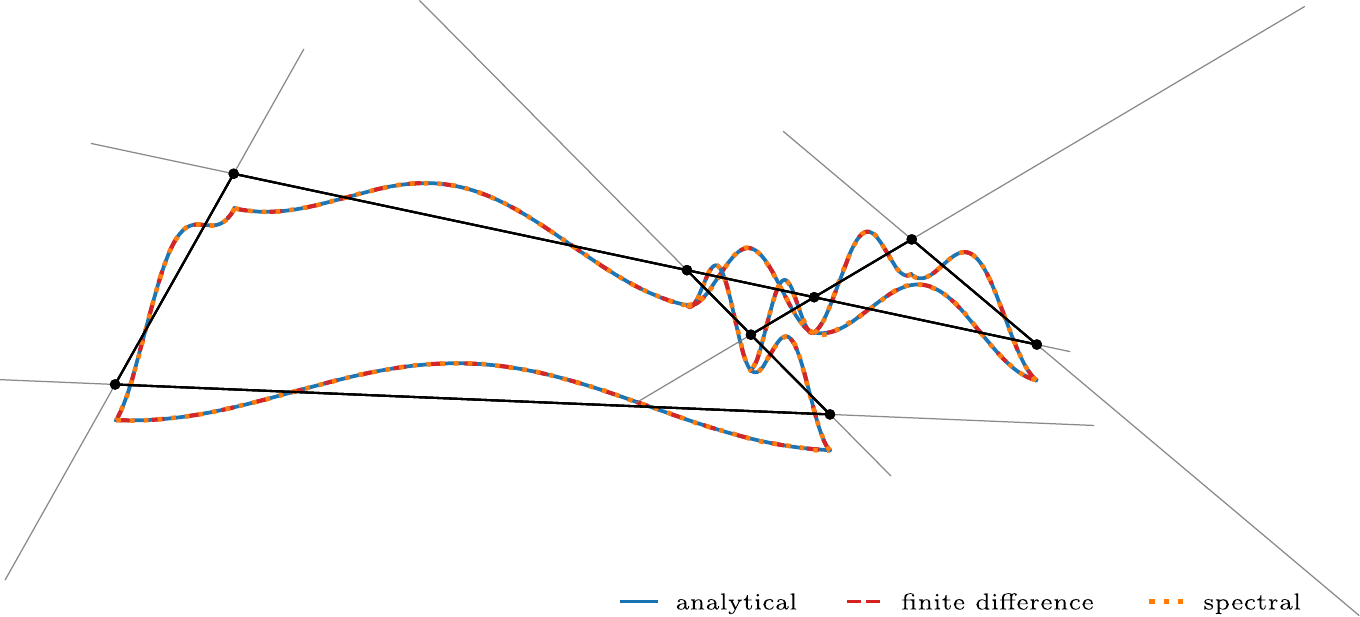}
    \caption{Solution of the Poisson equation on a metric random-line network with $N = 8$ nodes (black dots) and $M = 10$ edges (black lines). In the depicted realization of a random-line network, we independently position six ``needles'' (gray line segments) of unit length in the unit square. (See \cite{bottcher2020random} for further details about how to generate random-line networks.) 
    All edges have Kirchhoff boundaries. The source term on the right-hand side of the Poisson equation \eqref{eq:poisson} is $\rho_i(x_i) = \cos(2\pi x_i/\ell_i)/\ell_i^2$ for each edge. The solution of the Poisson equation is $\phi_i(x_i) = -\cos(2\pi x_i/\ell_i)/(4\pi^2)$ for each edge (solid blue curves). The dashed red and dotted orange curves indicate numerical solutions using finite-difference and spectral methods, respectively.}
    \label{fig:poisson_rlg}
\end{figure}

As a second metric network, we consider a hexagonal lattice with $N = 154$ nodes and $M = 213$ edges [see Fig.~\ref{fig:poisson_hexagonal}]. We set the lengths of all edges to $1$ and use $N_i = 1000$ discretization intervals for all edges $i \in \{1, \ldots, 213\}$ in the finite-difference approach. The source term is given by Eq.~\eqref{eq:source_term}. See our code repository~\cite{boettcher_gitlab} for the software implementation details for this example and all of our subsequent numerical examples.
We again observe the numerical simulations from both the finite-difference and spectral approaches closely resemble the analytical solution. Although we consider more than 200 edges and 1000 finite-difference discretizations per edge, using sparse-matrix solvers allows us to efficiently compute numerical solutions. In Appendix~\ref{app:larger_scale}, we consider metric hexagonal-lattice networks with up to about $10^4$ nodes and about $10^4$ edges.

As a third example of a metric network, we examine a random-line network with $N = 8$ nodes and $M = 10$ edges [see Fig.~\ref{fig:poisson_rlg}]. In a random-line network, one independently places line segments (\ie, ``needles'') of a specified length in a unit square (or other domain), forming an overlapping pattern~\cite{bottcher2020random}. Random-line networks, which are reminiscent of the Buffon needle graphs that were considered in prior works on metric networks~\cite{gaio2019nanophotonic,brio2022spectral}, are a useful toy model to studying PDEs on metric networks, as their edges are line segments of a specific length. Random-line networks and related spatial networks are relevant to the study of granular and particulate systems~\cite{papadopoulos2018network,nauer2020random,zorn2021charge}. 
%

\subsection{Heat equation}
\label{sec:heat}

\begin{figure}
    \centering
    \includegraphics{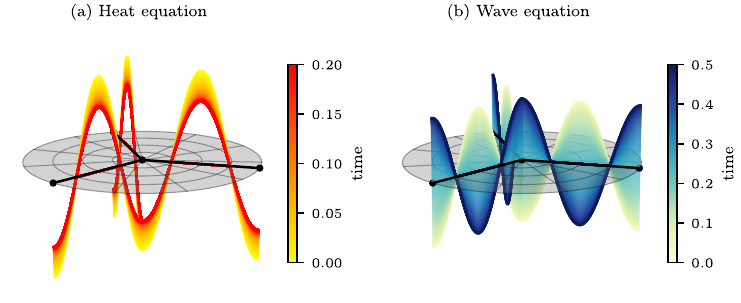}
    \caption{Solution of the heat and wave equations on a metric star network with $N = 4$ nodes (black dots) and $M = 3$ edges (black lines). All edges $i \in \{1,2,3\}$ have length $\ell_i = 1$ and Kirchhoff boundaries. (a) We consider the heat equation \eqref{eq:heat} with a source term $\rho_i(x_i) = \cos(2\pi x_i)$ and initial condition $\phi_i(x_i,0) = -3\cos(2\pi x_i)/(8\pi^2)$. In the limit as time $t \rightarrow \infty$, we recover the solution of the corresponding Poisson equation (\ie, $\lim_{t\rightarrow\infty} \phi_i(x_i,t) = -\cos(2\pi x_i)/(4\pi^2)$). (b) We consider the wave equation \eqref{eq:wave} with initial conditions $\phi(x_i,0) = -\cos(2\pi x_i)/(4\pi^2)$ and $\dot{\phi}(x_i,0) = 0$. The solution is $\phi(x_i,t) = -\cos(2\pi x_i)/(4\pi^2)\cos(2\pi t)$. The depicted solutions use a finite-difference approach with $N_i = 1000$ discretizations per edge.}
    \label{fig:heat_wave_star}
\end{figure}

We now consider the inhomogenous heat equation
\begin{equation}     \label{eq:heat}
    \frac{\partial \phi}{\partial t} = \tilde{\Delta} \phi - \rho
\end{equation}
with source term $\rho$. Its steady-state solutions satisfy the Poisson equation~\eqref{eq:poisson} in the limit as time $t \rightarrow \infty$.

In Fig.~\ref{fig:heat_wave_star}(a), we show the solution of the heat equation \eqref{eq:heat} for our 3-edge star network. We set $\rho_i(x_i) = \cos(2\pi x_i)$ and $\phi_i(x_i,0) = -3\cos(2\pi x_i)/(8\pi^2)$ for edges $i \in \{1,2,3\}$. For this initial condition and source term, $\lim_{t\rightarrow\infty} \phi_i(x_i,t) = -\cos(2\pi x_i)/(4\pi^2)$. In the Supplementary Information (SI)~\cite{animation1}, we show an animation of the evolution of $\phi_i(x,t)$ in Eq.~\eqref{eq:heat} for both our 3-edge star network and the hexagonal lattice from Section~\ref{sec:poisson}.
%


\subsection{Wave equation}
\label{sec:wave}
As a final example of a PDE on a metric network, we study the wave equation
\begin{equation}     \label{eq:wave}
    \frac{\partial^2 \phi}{\partial t^2}=\tilde{\Delta} \phi\,.
\end{equation}

In Fig.~\ref{fig:heat_wave_star}(b), we show the solution of the wave equation \eqref{eq:wave} for our 3-edge star network. For each edge $i \in\{1,2,3\}$, we set $\phi(x_i,0) = -\cos(2\pi x_i)/(4\pi^2)$ and $\dot{\phi}(x_i,0) = 0$. The solution is $\phi(x_i,t) = -\cos(2\pi x_i)/(4\pi^2)\cos(2\pi t)$. In the SI~\cite{animation2}, we show an animation of the evolution of $\phi_i(x,t)$ in Eq.~\eqref{eq:wave} for both our 3-edge star network and the hexagonal lattice from Section~\ref{sec:poisson}.
%


\section{Conclusions and discussion}
\label{sec:conclusion}
Metric networks give a mathematical framework to study spatially extended dynamics, such as partial differential equations, on networked systems.
Metric networks have applications in a variety of scientific fields, including modeling the mechanical properties of materials~\cite{hrennikoff1941solution,ashurst1976microscopic,beale1988elastic,hassold1989brittle,herrmann1989fracture,gusev2004finite,bottcher2021computational}, describing quantum dynamics in thin structures~\cite{ruedenberg1953free,alexander1983superconductivity,kottos1997quantum,kottos1999periodic,kuchment2002graph,kuchment2003quantum,kuchment2005quantum,berkolaiko2013introduction,Berkolaiko2017}, analyzing information propagation in transmission lines~\cite{paul2007analysis,strub2019modeling,alonso2017power,chen2017power,muranova2019notion,muranova2020eigenvalues,muranova2021effective,muranova2022effective,muranova2022networks}, and simulating gas flow in pipelines~\cite{banda2006gas,brouwer2011gas,domschke2011combination,mindt2019entropy,DomschkeHillerLangetal.2021}. 

Over the last two to three decades, the study of metric networks has progressed in parallel to (and largely independent of) developments in more conventional network science. In network science, analysis of the interplay between network structure and dynamics has long been a key topic. Such research has used combinatorial networks, rather than metric networks, and accordingly it has focused primarily on ordinary differential equations and (both deterministic and stochastic) agent-based models on networks, rather than on PDEs.
Including intervals on edges with a metric structure yields a natural setting to analyze PDEs on networked systems.
Unlike in ODE dynamics, dynamical processes on metric networks require one to specify continuity conditions and boundary conditions in addition to an initial condition. The interplay between network structure and boundary conditions is associated with different possible characteristic eigenmodes of a metric network.

In the present paper, we overviewed several analytical and numerical approaches that are essential for studying fundamental linear PDEs, such as the time-independent Schr\"odinger equation, the Poisson equation, the heat equation, and the wave equation. We expanded on the spectral approach of \cite{gaio2019nanophotonic,brio2022spectral} to account for degenerate eigenmodes and various algorithm inputs, including the range of characteristic wavenumbers, the optimization bounds, and the rounding precision of potentially equivalent characteristic wavenumbers. Although a spectral approach is useful to identify eigenmodes in a metric network, it may be challenging to accurately determine a large number of (potentially degenerate) eigenmodes to unambiguously determine the solution of a given PDE problem.

Complementing the numerical results in \cite{brio2022spectral}, which focused on problems involving the Poisson equation and the telegraph equation on a metric network with three nodes,\footnote{They obtained solutions of PDEs on metric networks with 3 nodes, and they computed wavenumbers for metric networks with up to about 100 nodes.} we examine the Poisson equation, heat equation, and wave equation on three distinct metric networks with larger numbers of nodes and edges. The spectral solver, finite-difference solvers, and visualization routines that we developed in our investigation are available at \url{https://gitlab.com/ComputationalScience/metric-networks}.

There are numerous worthwhile research directions to pursue in future work. Given the challenges of obtaining accurate solutions of different types of PDEs on metric networks, it is important to further develop and improve numerical solvers. In addition to numerical techniques like finite-difference, finite-element, finite-volume, and spectral methods, potential approaches can also encompass physics-informed neural networks (PINNs) \cite{blechschmidt2022comparison}, including ones that use spectral information~\cite{xia2023spectrally}. In the spirit of \cite{golubitsky2023dynamics}, another promising avenue is extending
symmetry arguments~\cite{band2017quotients,jevzek2021application} to various families of metric networks.
 The analysis of symmetries can help improve understanding of how specific structural features influence eigenmodes and their degeneracies. It may also be worthwhile to explore the connections between studies of metric networks and related research on problems such as the topological Dirac equation on networks and simplicial complexes~\cite{bianconi2021topological}.
It is also important to study nonlinear PDEs on metric networks.

%


\section*{Acknowledgments}
We thank Pia Domschke and Hannah Kravitz for helpful comments. MAP thanks Leonid Bunimovich for introducing him to quantum graphs more than 20 years ago. It has taken a long time, but that introduction created the earliest kernel that has finally led to the present paper.
%


\section*{Code availability}
Our code is publicly available at \url{https://gitlab.com/ComputationalScience/metric-networks}.


\clearpage
\appendix

\section{Symmetries}
\label{app:symmetries}
It is very important to consider symmetries to understand dynamical processes on networks~\cite{golubitsky2023dynamics}. 
In this section, we will illustrate how to use symmetry groups~\cite{weyl1950theory,dresselhaus2007group,zee2016group,Bouchard_2020} to identify degeneracies in the eigenmodes of a metric network. Such analytical insights can help one assess whether or not a numerical method has successfully identified
all eigenmodes. The essential idea is to determine the irreducible representations of a symmetry group and then compute eigenmode degeneracies by considering the dimensions of the relevant irreducible representations.

Degenerate eigenmodes, such as the ones in Eq.~\eqref{eq:star_degenerate}, are usually associated with a symmetry of a metric network~\cite{berkolaiko2017simplicity}. For example, in our 3-edge star network (see Fig.~\ref{fig:star_graph}), we obtain the same characteristic wavenumbers and corresponding eigenmodes if we permute the three identical length-$\ell$ edges, which have the same function space and operator space. 
In this example, the relevant symmetry group is the symmetric group $S_3$, which consists of the $|S_3| = 3! = 6$ possible permutations of the elements of the set $\{1,2,3\}$. Using cycle notation,\footnote{In cycle notation, one describes a permutation as a product of disjoint cycles~\cite{zee2016group}. 
In each cycle, one rearranges a set of elements among themselves. For example, in cycle notation, the permutation $(1,2,3)$ indicates that we map $1$ to $2$, $2$ to $3$, and $3$ to $1$.} the set of elements of $S_3$ is $\{\mathds{1}, (1,2), (2,3), (1,3), (1,2,3), (1,3,2)\}$, where $\mathds{1}$ is the identity element. Permutations that involve two elements are called ``transpositions", and permutations that involve three elements are called ``3-cycles". The three distinct types of cycle structures yield $C = 3$ conjugacy classes, which are relevant for characterizing degenerate eigenmodes.

Consider the permutation representation $P \colon S_3\rightarrow \mathrm{GL}(\mathbb{R}^3)$ of $S_3$ in which the edges $e_1$, $e_2$, and $e_3$ are represented by the vectors
\begin{equation}
    e_1=
    \begin{pmatrix}
        1\\
        0\\
        0
    \end{pmatrix}\,,\quad
    e_2=
    \begin{pmatrix}
        0\\
        1\\
        0
    \end{pmatrix}\,,\quad 
    e_3=
    \begin{pmatrix}
        0\\
        0\\
        1
    \end{pmatrix}\,.
\end{equation}
In this representation, the permutation $\pi = (1,2,3)$ is 
\begin{equation}
    P(\pi)=
    \begin{pmatrix}
        0 & 0 & 1 \\
        1 & 0 & 0 \\
        0 & 1 & 0
    \end{pmatrix}\,.
\end{equation}
Observe that $P(\pi)e_1 = e_2$, $P(\pi)e_2 = e_3$, and $P(\pi)e_3 = e_1$. One can similarly determine matrix representations of the remaining five elements of the group $S_3$. 

Given a representation $P$, the character $\chi^{(P)}\colon S_3\rightarrow\mathbb{R}$ assigns the trace of the corresponding matrix representation to each group element $g \in S_3$. That is,
\begin{equation}
    \chi^{(P)}(g)=\mathrm{tr}(P(g))\,.
\end{equation}
For the three conjugacy classes of $S_3$ and the permutation representation $P$, the corresponding characters are $\chi_1^{(P)} = 3$ (the identity permutation, in which no elements are rearranged), $\chi_2^{(P)} = 1$ (the transpositions, in which two elements are swapped and one element remains in its current position), and $\chi_3^{(P)} = 0$ (the 3-cycles, in which no elements remain in their current positions). 

A representation $P$ is ``semisimple'' (\ie, completely reducible) if one can decompose it into a direct sum of irreducible representations $P^{(\alpha)}$. That is,
\begin{equation}     \label{eq:decomposition}
    P = \bigoplus_{\alpha = 1}^C c_\alpha P^{(\alpha)}\,,
\end{equation}
where $C$ is the number of conjugacy classes and $c_\alpha$ denotes a decomposition coefficient. Taking the trace of Eq.~\eqref{eq:decomposition} yields
\begin{equation}
    \chi^{(P)}(g) = \sum_{\alpha = 1}^C c_\alpha \chi^{(\alpha)}(g)\,.
\end{equation}

The dimensions $d_{\alpha}$ of the irreducible representations $P^{(\alpha)}$ of a finite group $G$ satisfy
\begin{equation}
    \sum_{\alpha=1}^C d_\alpha^2 = |G|\,.
\end{equation}
For the group $S_3$, the only dimensions of the three irreducible representations that satisfy $d_1^2 + d_2^2 + d_3^3 = 6$ are $d_1 = d_2 = 1$ and $d_3 = 2$. The two one-dimensional irreducible representations correspond to the trivial and sign representations. In the trivial representation, one maps each element of $S_3$ to $1$. In the sign representation, one maps each permutation to its corresponding sign, which is $1$ for even permutations and $-1$ for odd permutations. We will show below that the two-dimensional irreducible representation leads to the eigenmode degeneracy that we observed in Section~\ref{sec:star_graph}.
This representation is the ``standard representation'' of $S_3$. 
It is a ``faithful" representation, which means that it gives a one-to-one mapping of group elements to their corresponding matrices. By contrast, the other two representations are not faithful.

A representation $P$ of a finite group $G$ is irreducible if and only if
\begin{equation}
    \sum_{\alpha=1}^{C} n_\alpha |\chi_\alpha^{(P)}|^2 = |G|\,,
\end{equation}
where $n_\alpha$ denotes the number of elements in the $\alpha$th conjugacy class. For a reducible representation, 
\begin{equation}
    \sum_{\alpha=1}^{C} n_\alpha |\chi_\alpha^{(P)}|^2 > |G|\,.    
\end{equation}

Inserting the values of $n_\alpha$ and $\chi_\alpha^{(P)}$ that are associated with the permutation representation yields
\begin{equation}
    \sum_{\alpha = 1}^{C} n_\alpha |\chi_\alpha^{(P)}|^2 = 1 \times 3^2 + 3 \times 1^2+2\times 0^2 = 12\,.    
\end{equation}
Because $12 > |S_3| = 6$, the permutation representation is reducible. The decomposition coefficients $c_\alpha$ [see Eq.~\eqref{eq:decomposition}] are 
\begin{equation}     \label{eq:decomposition_coefficients}
    c_\alpha = \frac{1}{|G|}\sum_{\beta = 1}^C n_{\beta} \chi_\beta^{(P)} (\chi_\beta^{(\alpha)})^\ast\,.
\end{equation}

\begin{table}
\centering
\begin{tabular}{l | c c c}
\toprule
 & $\mathds{1}$ & $(1,2)$ & $(1,2,3)$ \\ \hline
trivial representation & 1 & 1 & 1 \\ 
sign representation & 1 & $-1$ & 1 \\ 
standard representation & 2 & 0 & $-1$ \\
\bottomrule
\end{tabular}
\caption{Character table of the symmetric group $S_3$.}
\label{tab:S3_characters}
\end{table}

Using the character table (see Table \ref{tab:S3_characters}) yields
\begin{align}
    c_1 &= \frac{1}{6}\left(1\times 3 \times 1 +3\times1\times 1 + 2\times 0\times 1\right) = 1 \,, \\
    c_2 &= \frac{1}{6}\left(1\times 3 \times 1+3\times 1 \times (-1) + 2\times 0\times 1\right) = 0 \,, \\
    c_3 &= \frac{1}{6}\left(1\times 3 \times 2 + 3\times 1 \times 0 + 2 \times 0 \times (-1)\right) = 1\,.
\end{align}
We can thus decompose the permutation representation $P$ of the symmetric group $S_3$ as
\begin{equation}
    P = P^{(1)}\oplus P^{(3)}\,.
\end{equation}
That is, the permutation representation $P$ of the symmetric group $S_3$ is the direct sum of the trivial irreducible representation $P^{(1)}$ and the two-dimensional irreducible representation $P^{(3)}$.

For the linear time-independent Schr\"odinger equation on the 3-edge metric star network (see Section~\ref{sec:star_graph}), the characteristic wavenumbers and eigenmodes that are associated with the Hamiltonian $\widetilde{\mathcal{H}} = -\tilde{\Delta}$ (\ie, the generalized negative Laplacian that includes continuity and boundary conditions) are invariant with respect to permutations of the edges.
That is, the permutation operator $P$ commutes with $\widetilde{\mathcal{H}}$ (\ie, $[\widetilde{\mathcal{H}},P] = 0$), so
\begin{equation}
    \widetilde{\mathcal{H}}P f=P\widetilde{\mathcal{H}}f = -Pk^2f =-k^2 Pf\,.
\end{equation}
We now choose a basis so that the permutation operator decomposes into a direct sum of the irreducible representations $P^{(1)}$ and $P^{(3)}$. By Schur's Lemma, the Hamiltonian $\widetilde{\mathcal{H}}$ becomes diagonal in this basis.
\begin{lemma}[Schur's Lemma]
If $P \colon G\rightarrow \mathrm{GL}(V)$ is an irreducible representation of a finite group $G$ and there exists a matrix $E$ that commutes with every element $g \in G$ (\ie, $E P(g) = P(g)E$ for all $g\in G$), then $E = \lambda I$, where $I$ is the identity matrix and $\lambda \in \mathbb{C}$.
\end{lemma}

According to Schur's lemma, in the basis in which the representation $P$ decomposes into irreducible representations $P^{(1)}$ and $P^{(3)}$, the basis vectors are eigenstates of the Hamiltonian $\widetilde{\mathcal{H}}$. The dimensions $d_1 = 1$ and $d_3 = 2$ of these irreducible representations correspond to the degeneracies that we observed in the star network with three length-$\ell$ edges~\cite{nussbaum1968group,zee2016group,Bouchard_2020}. 

We now apply the same reasoning to a metric star network with four length-$\ell$ edges. Without explicitly calculating the characteristic wavenumbers and eigenmodes, we deduce the underlying degeneracies by examining the irreducible representations that are associated with the permutation representation of $S_4$. 
A similar calculation as with $S_3$ shows that the permutation representation of $S_4$ is the direct sum of the one-dimensional trivial representation and the three-dimensional standard representation. Therefore, the eigenmode degeneracies are $1$ and $3$. 
Indeed, an explicit calculation of the determinant of the coupling-condition matrix $T_{\fourpointstar{0.15cm}}(k_m)$ [see Eq.~\eqref{eq:homogenous_system}] for this 4-edge star network yields $\mathrm{det}(T_{\fourpointstar{0.15cm}}(k_m)) = 4 \cos(k_m \ell)^3 \sin(k_m \ell)$,\footnote{The determinant of the coupling-condition matrix for an $M$-edge metric star network with edges $i \in \{1,\ldots,M\}$ of length $\ell_i$ is $\mathrm{det}(T_{\fourpointstar{0.15cm}}(k_m)) = \sum_{i = 1}^M\tan(k_m \ell_i) \prod_{i = 1}^M \cos(k_m \ell_i)$~\cite{swindle2019spectral}.} 
demonstrating that the group-theoretically determined degeneracies coincide with the ones that we determined by calculating the determinant of $T_{\fourpointstar{0.15cm}}(k_m)$.

In summary, when using a group-theoretic approach to determine eigenmode degeneracies for the Schr\"odinger equation on a metric network {(and, more generally, to determine the eigenstate degeneracies of a Hamiltonian that is associated with a metric network)},
we follow the following procedure:
\begin{itemize}
    \item Given a metric network, determine the characters of the relevant representation of its symmetry group. (A relevant representation describes the action of the symmetry group on the system under consideration.)
    \item Calculate the decomposition coefficients $c_\alpha$ [see Eq.~\eqref{eq:decomposition_coefficients}] and the corresponding decomposition into irreducible representations [see Eq.~\eqref{eq:decomposition}].
    The dimensions of the irreducible representations with nonzero coefficients $c_\alpha$ are equal to the degeneracies of the system's eigenstates.
\end{itemize}

It is possible for accidental symmetries to cause equality of eigenvalues that are associated with different irreducible representations~\cite{mcintosh1959accidental}. For a detailed treatment of symmetry operations on metric networks and their decomposition into substructures (so-called ``quotient graphs"), see \cite{band2017quotients,jevzek2021application}.

This caveat notwithstanding, degeneracies in the eigenvalues of a PDE on a metric network are usually associated with the symmetries of the metric network. One can use small perturbations of the edge lengths to break symmetries~\cite{Berkolaiko2017}. For Schr\"odinger equations on metric networks with Kirchhoff flux boundary conditions, the
eigenvalues are usually simple (\ie, their algebraic multiplicity is usually $1$)~\cite{friedlander2005genericity,colin2015semi}. 
Nevertheless, given the broad relevance of symmetric network structures, it is important to consider the connections between symmetry groups and eigenmode degeneracies.
%


\section{Large networks} \label{app:larger_scale}
\begin{figure}
    \centering
    \includegraphics{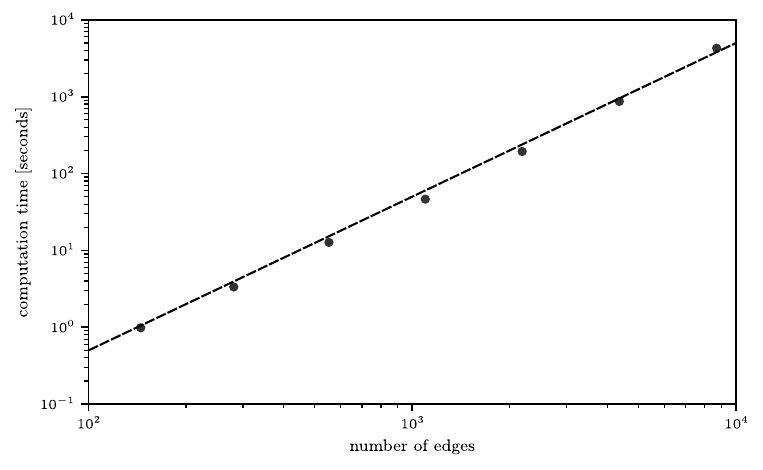}
    \caption{Computation time for solving the Poisson equation \eqref{eq:poisson} on metric hexagonal lattices with a finite-difference approach for different numbers $N$ of nodes and numbers $M$ of edges: $(N,M) \in \{(106,145),(202,281),(394,553),(778,1097),(1546,2185),(3082,4361),(6154,8713)\}$. 
        All edges $i \in \{1,\ldots,M\}$ have length $\ell_i = 1$ and Kirchhoff boundaries. The source term on the right-hand side of the Poisson equation \eqref{eq:poisson} is $\rho_i(x_i) = \cos(2\pi x_i)$ for each edge. The  solution of the Poisson equation is $\phi_i(x_i) = -\cos(2\pi x_i)/(4\pi^2)$ for each edge. The dashed black line corresponds to a power law with exponent $2$.
        }
    \label{fig:runtime_edges}
\end{figure}
We solve the Poisson equation \eqref{eq:poisson} on metric hexagonal-lattice networks with different numbers of nodes and edges. The source term is given by Eq.~\eqref{eq:source_term}. 
We set the lengths $\ell_i$ of all edges $i \in \{1,\ldots,M\}$ to $1$, and all boundaries are of Kirchhoff type. The solution of the corresponding Poisson equation is $\phi_i(x_i) = -\cos(2\pi x_i)/(4\pi^2)$. 

To solve the Poisson equation on these metric networks, we employ a finite-difference approach and set the number of discretization intervals to $N_i = 1000$ for all edges. In Fig.~\ref{fig:runtime_edges}, we show the computation time as a function of the number of edges. The largest metric network that we consider has $6154$ nodes and $8713$ edges. Solving the Poisson equation on this network takes about $4300$ seconds (\ie, about 1.2 hours) on one i7 CPU core with a 1.8 GHz clock speed. In all of our simulations, the mean-squared error between the numerical and analytical solutions is less than $10^{-17}$. 
For this error calculation, we use vectors that contain discretized solutions of the Poisson equation at all edges.


%
\bibliographystyle{siamplain}
\bibliography{references}
%


\end{document}